\documentclass[11pt]{article}

\usepackage[T1]{fontenc}
\usepackage[latin1]{inputenc}
\usepackage{amsmath,amsthm,amssymb}
\usepackage{mathabx}

\usepackage[margin=1in]{geometry}
\usepackage{xcolor}
\usepackage{hyperref}

\newtheorem{prop}{Proposition}[section]
\newtheorem{cor}[prop]{Corollary}
\newtheorem{lem}[prop]{Lemma}
\newtheorem{theo}[prop]{Theorem}

\def\ricc{\mbox{\rm Ricc}}
\def\tr{\mbox{\rm tr}}

\newcommand{\EE}{\mathbb{E}}
\newcommand{\LL}{\mathbb{L}}
\newcommand{\RR}{\mathbb{R}}
\newcommand{\SB}{\mathbb{S}}

\newcommand{\Ca}{ {\cal C }}
\newcommand{\Ea}{ {\cal E }}
\newcommand{\Oa}{ {\cal O }}
\newcommand{\Ra}{ {\cal R }}
\newcommand{\Fa}{ {\cal F }}
\newcommand{\Pa}{ {\cal P }}

\begin{document}

% =======================================================
\title{On the Stability of Kalman-Bucy Diffusion Processes\thanks{Please cite this work: ``A.N. Bishop and P. Del Moral. On the Stability of Kalman-Bucy Diffusion Processes. SIAM Journal on Control and Optimization. 55(6):4015--4047 (2017); arxiv e-print arxiv.org/abs/1610.04686 updated.''}}
\author{Adrian N. Bishop and Pierre Del Moral}
\date{}

\maketitle

\begin{abstract}
The Kalman-Bucy filter is the optimal state estimator for an Ornstein-Uhlenbeck diffusion given that the system is partially observed via a linear diffusion-type (noisy) sensor. Under Gaussian assumptions, it provides a finite-dimensional exact implementation of the optimal Bayes filter. It is generally the only such finite-dimensional exact instance of the Bayes filter for continuous state-space models. Consequently, this filter has been studied extensively in the literature since the seminal 1961 paper of Kalman and Bucy. The purpose of this work is to review, re-prove and refine existing results concerning the dynamical properties of the Kalman-Bucy filter so far as they pertain to filter stability and convergence. The associated differential matrix Riccati equation is a focal point of this study with a number of bounds, convergence, and eigenvalue inequalities rigorously proven. New results are also given in the form of exponential and comparison inequalities for both the filter and the Riccati flow. \medskip

\emph{Keywords}: 
differential Riccati equations, diffusion flows, Kalman-Bucy diffusion, Kalman-Bucy filter, transition semigroups.\medskip

\end{abstract}

\setcounter{tocdepth}{2}
{\small \tableofcontents}
% =======================================================

\section{Introduction}

The aim of this study is to review, reprove, and also to refine a number of existing stability results on the Kalman-Bucy filter and the associated Riccati equation. We correct prior work where necessary. New results are also given in the form of exponential and comparison inequalities for both the stochastic flow of the filter, and the Riccati flow. \textbf{This work is intended to be a complete and self-contained analysis on the stability and convergence of Kalman-Bucy filtering}; with detailed proofs of each necessary result.

Consider a linear-Gaussian filtering model of the following form
\begin{equation}\label{lin-Gaussian-diffusion-filtering}
\left\{
\begin{array}{rcl}
	dX_t &=&A_t\,X_t~dt~+~R^{1/2}_{1}~dW_t\\
	dY_t &=&C_t\,X_t~dt~+~R^{1/2}_{2}~dV_{t}.
\end{array}
\right.
\end{equation}
Here, $(W_t,V_t)$ is an $(r_1+r_2)$-dimensional standard Brownian motion. Let $\Fa_t=\sigma\left(Y_s,~s\leq t\right)$ be the $\sigma$-algebra filtration generated by the observations and $Y_0=0$. Assume $X_0$ is a $r_1$-valued independent random vector with mean $\EE(X_0)$ and finite covariance $P_0$. Note $X_0$ is not necessarily Gaussian.

Further, $A_t$ is a square $(r_1\times r_1)$-matrix, $C_t$ is an $(r_2\times r_1)$-matrix, and $R^{1/2}_{1}$ and $R^{1/2}_{2}$ are symmetric $(r_1\times r_1)$ and $(r_2\times r_2)$ matrices. The eigenvalues of $A,C,R_{1},R_{2}$ are bounded above and below (uniformly in time) and those of $R_{1},R_{2}$ are uniformly bounded positive.

We consider both time-varying (e.g. $A_t$) and time-invariant signal models (e.g. $A_t=A$) and the convergence properties of the respective filters and associated Riccati equations. Typically, we state general results for the time-varying signal first, and follow this with more quantitative results in the time-invariant case.

When $X_0$ is Gaussian, it is well-known that the conditional distribution of the signal state $X_t$ given $\Fa_t$ is a $r_1$-dimensional Gaussian  distribution with a mean and covariance matrix 
$$
\widehat{X}_t:=\EE(X_t~|~\Fa_t)\quad\mbox{\rm and}\quad
P_t:=\EE\left(\left(X_t-\EE(X_t~|~\Fa_t)\right)\left(X_t-\EE(X_t~|~\Fa_t)\right)^{\prime}\right)
$$ 
given by the Kalman-Bucy and the Riccati equations 
\begin{eqnarray}
d\widehat{X}_t&=&A_t~\widehat{X}_t~dt+P_{t}~C_t^{\prime}R^{-1}_{2}~\left(dY_t-C_t\widehat{X}_tdt\right)\label{nonlinear-KB-mean}\\
\partial_tP_t&=&\ricc(P_t) \label{nonlinear-KB-Riccati}
\end{eqnarray}
with the Riccati drift function defined for any positive definite $Q$ by
$$
\ricc(Q)=A_tQ+QA_t^{\prime}-QS_tQ+R_1\quad\mbox{\rm with}\quad S_t:=C_t^{\prime}R^{-1}_{2}C_t.
$$
Note that $S_t$ is positive semi-definite and time-varying whenever, e.g., $C_t$ is time-varying. We may take $Q$ only positive semi-definite, but for simplicity throughout we assume $Q$ positive definite.

The Kalman-Bucy filter is the $\LL_2$-optimal state estimator for an Ornstein-Uhlenbeck diffusion given that the system is partially observed via a linear diffusion-type (noisy) sensor; see \cite{kalman61,bucy68}.

Stability of this filter was initially studied by R.E. Kalman and R.S. Bucy in their seminal paper \cite{kalman61}, with related prior work by Kalman \cite{kalman60,kalman60-2,kalman60-3} and later work by Bucy \cite{bucy2}. In \cite{anderson71} the stability of this filter was analysed under a relaxed controllability condition. It was analysed again in \cite{ocone-pardoux} for systems with non-Gaussian initial state via Kallianpur-Striebel-type change of probability measures. An alternative approach is to  consider the following conditional nonlinear McKean-Vlasov-type diffusion process
\begin{equation}\label{Kalman-Bucy-filter-nonlinear-ref}
d\overline{X}_t=A_t~\overline{X}_t~dt~+~R^{1/2}_{1}~d\overline{W}_t+\Pa_{\eta_t}C_t^{\prime}R^{-1}_{2}~\left[dY_t-\left(C_t\overline{X}_tdt+R^{1/2}_{2}~d\overline{V}_{t}\right)\right]
\end{equation}
where  $(\overline{W}_t,\overline{V}_t,\overline{X}_0)$ are independent copies of $(W_t,V_t,X_0)$ (thus independent of
 the signal and the observation path). In the above displayed formula $\Pa_{\eta_t}$ stands for the covariance matrix
$$
\Pa_{\eta_t}=\eta_t\left[(e-\eta_t(e))(e-\eta_t(e))^{\prime}\right]
\quad\mbox{\rm with}\quad \eta_t:=\mbox{\rm Law}(\overline{X}_t~|~\Fa_t)\quad\mbox{\rm and}\quad
e(x):=x.
$$
We shall call this probabilistic model the Kalman-Bucy (nonlinear) diffusion process. 

The nonlinear interaction does not take place only on the drift part, but also on the diffusion matrix functional. In addition the nonlinearity does not depend on the distribution of the random states $\pi_t=\mbox{\rm Law}(\overline{X}_t)$ but on their conditional distributions $\eta_t:=\mbox{\rm Law}(\overline{X}_t~|~\Fa_t)$. The well-posedness of this nonlinear diffusion is discussed in~\cite{dm-16-enkf}. 

The nonlinear Kalman-Bucy diffusion is, in some sense, a generalized description of the Kalman-Bucy filter. The Riccati equation (\ref{nonlinear-KB-Riccati}) is captured in the nonlinear term of the diffusion. 
More precisely, the conditional expectations of the random states $\overline{X}_t$ and their conditional covariance matrices $\Pa_{\eta_t}$ w.r.t. $\Fa_t$ satisfy the Kalman-Bucy and the Riccati equations \eqref{nonlinear-KB-mean} and \eqref{nonlinear-KB-Riccati}, {\em even when the initial state is not Gaussian}. That is, if we {\em redefine} 
 \begin{equation}\label{ce-intro}
	\widehat{X}_t := \EE\left(\overline{X}_t~|~\Fa_t\right) \qquad\mbox{\rm and} \qquad P_t := \Pa_{\eta_t}
\end{equation} 
then the flow of this (conditional) mean and covariance satisfy \eqref{nonlinear-KB-mean} and \eqref{nonlinear-KB-Riccati} irregardless of the distribution of $X_0$; see \cite{dm-16-enkf}. We assume this more general definition of $\widehat{X}_t$ and $P_t$ when referring to \eqref{nonlinear-KB-mean} and \eqref{nonlinear-KB-Riccati} going forward. 

Note that the flow of matrices $\Pa_{\eta_t}$ depends only on the covariance matrix of the initial state $\overline{X}_0$. This property follows from the specially designed structure of the nonlinear diffusion, which ensures that the mean and covariance matrices satisfy the Kalman-Bucy filter and Riccati equations. This structure simplifies the stability analysis of this diffusion. Given $\Pa_{\eta_0}$ the Kalman-Bucy diffusion \eqref{Kalman-Bucy-filter-nonlinear-ref} can be interpreted as a non-homogeneous Ornstein-Uhlenbeck type diffusion with a conditional covariance matrix $P_t=\Pa_{\eta_t}$ that satisfies the Riccati equation \eqref{nonlinear-KB-Riccati} starting from $P_0=\Pa_{\eta_0}$. In this sense, the nonlinearity of the process is encapsulated by the Riccati equation.

Analysis of this diffusion allows one to capture non-Gaussian initial states even for time-varying signal models. This class of nonlinear diffusion also arises in the mathematical and the numerical
foundations of ensemble-Kalman-Bucy filters and data assimilation~\cite{evensen-review}. In this context, the stability properties of the Kalman-Bucy diffusion are essential for analyzing the long-time behaviour of this class of algorithm.

Reiterating, in this work, we revisit the stability of the Kalman-Bucy filter and we study for the first time the stability properties of the Kalman-Bucy diffusion. We derive new exponential inequalities detailing the convergence of the filter and the diffusion with arbitrary initial conditions, and the convergence properties of the associated differential Riccati equation. The classical study of Riccati equations in control and estimation theory is motivated by their relationship with Kalman-Bucy filtering and linear-quadratic optimal control theory \cite{kalman61,kalman60,kalman60-2}. Indeed, the two topics are dual, and the two relevant differential Riccati equations are (mostly) equivalent up to a time-reversal. We deal here primarily with the forward-type equation associated with the evolution of the Kalman-Bucy filtering error.

We review now some of the key literature on the (deterministic) matrix Riccati-type differential equation, i.e. quadratic matrix differential equations \cite{kalman61}. Our interest in this equation follows because it describes the covariance flow of the Kalman-Bucy state estimation error. However, the properties and behaviour of this equation are of interest in their own rights. Bucy \cite{bucy2} originally studied a number of global properties of the differential matrix Riccati equation. In particular, he proved that solutions exist for all time when the initial condition is positive semi-definite, he proved a number of important monotonicity properties, along with bounds on the solution stated in terms of the controllability and observability Gramians. Bucy \cite{bucy2} studied when the solution of the autonomous Riccati equation converges to a solution of an associated (fixed-point) algebraic Riccati equation, and finally he proved exponential stability of the time-varying Kalman-Bucy filter along with an exponential forgetting property of the associated Riccati equation. We review and re-prove these results here via novel methods. We also refine quantitative estimates.

It is worth noting some history concerning Bucy's uniform bounds. The original upper and lower bounds on the (time-varying) Riccati equation given in \cite{bucy2} were particularly elegant in appearance; being given in terms of the relevant observability and controllability Gramians. However, as noted in \cite{hitz72}, there was a crucial (yet commonly made) error in the proof which invalidated the result as given. This error was repeated (and/or overlooked) in numerous subsequent works; including by the current authors in the first writing of this work. A correction \cite{bucy72corrected} was noted in a reply to \cite{hitz72}; see Bucy's reply \cite{bucy72remarks} and a separate reply by Kalman \cite{kalman72}. However, a complete reworking of the result did not appear in entirety, it seems, until much later in \cite{delyon2001}. We remark that, in some sense, the qualitative nature of the Kalman-Bucy filter's stability was not jeopardised; as noted by Bucy \cite{bucy72remarks,bucy72corrected} and Kalman \cite{kalman72}. However, given time-varying signal models, the lack of a complete proof on the uniform boundedness of the Riccati equation in quantitative terms was somewhat unsatisfactory.

Associated with the differential Riccati equation is the (fixed-point) algebraic Riccati equation whose solution(s) correspond to the equilibrium point(s) of the corresponding differential equation. This algebraic equation was studied by Bucy in \cite{bucy72} and it was shown that there exists an unstable negative definite solution (in addition to the desired positive definite equilibrium). A detailed study of the algebraic Riccati equation was given by J.C. Willems \cite{willems71} who considered characterising all solutions. Bucy \cite{bucy75} later considered the so-called structural stability of these solutions. Detectability and stabilisability conditions are necessary and sufficient for a unique stabilising positive semi-definite solution of the algebraic equation \cite{kucera72}. See also \cite{kucera72,callier81,wimmer85} and the early review paper \cite{kucera73} for related literature. A (marginally) stable solution of the algebraic Riccati equation exists under detectability conditions; see \cite{Molinari73,Molinari77,wimmer85,Poubelle86}. The discussion in \cite[Chapter 2 and 3]{Bittanti91} is also of general interest here; as is \cite{Lancaster1995}.

Returning to the differential Riccati equation \cite{bucy2}, convergence to a stabilising fixed-point was studied extensively in \cite{callier81}, where a number of generalised convergence conditions in terms of the system model were given. We also note the early paper \cite{wonham68} that studied convergence and dealt further with a generalized version of the Riccati equation with a linear perturbation term. See also the seminal text \cite{Kwakernaak72}. A geometric analysis of the differential Riccati equation and its solution(s) is given in \cite{shayman86}. In \cite{gevers85} sufficient conditions are given such that the solution of the differential Riccati equation at any instant is stabilising; see also \cite{Poubelle88}. In other work \cite{Poubelle86,nicolao92,callier95,Park97} convergence to a (marginally) stable solution was studied again, with further relaxations and where necessary conditions on the system model were addressed.

Finally, we point to the texts \cite{reid72,Bittanti91,abou-kandil03}, dedicated to the Riccati equation, for further background and results (many of which are tangent to the discussion relevant here).

Given convergence of the differential Riccati flow and some associated semigroups, one typically concludes, in a straightforward way, the corresponding stability of the Kalman-Bucy filter; see the work of Bucy \cite{kalman61,bucy2} and the studies \cite{anderson71,ocone-pardoux}. However, we refine this conclusion in this work, with exponential inequalities and some related results.

\subsection{Statement of the main results and paper organization}

Let $\left\Vert\cdot\right\Vert_2$ be the Euclidean norm on $\RR^{r}$, or the spectral norm on $\RR^{r\times r}$,
for some $r\geq 1$. We denote by $\SB_r$ the set of $(r\times r)$ real symmetric matrices, and by $\SB_r^+$ the subset of positive-definite matrices. To describe our main results with some precision we need to introduce some notation. For any $0\leq s\leq t$, 
and $(x,Q)\in (\RR^{r_1}\times \SB^+_{r_1})
$ we let 
$$
\left(\varphi_{s,t}(x),\psi_{s,t}(x,Q),\overline{\psi}_{s,t}(x,Q),\phi_{s,t}(Q)\right)\in (\RR^{r_1}\times\RR^{r_1}\times\RR^{r_1}\times \SB^+_{r_1})
$$
 be respectively the  flow of signal (\ref{lin-Gaussian-diffusion-filtering}), the Kalman-Bucy filter (\ref{nonlinear-KB-mean}), the Kalman-Bucy diffusion (\ref{Kalman-Bucy-filter-nonlinear-ref}), and the Riccati equation (\ref{nonlinear-KB-Riccati}).  We take here the conventional observability/controllability conditions as holding; see the Standing Assumption (\ref{steady-state-eq-2}) in Section \ref{desc-sec-intro}.
 
 In Section \ref{desc-sec-intro} we introduce the relevant signal model, the Kalman-Bucy filter and an associated nonlinear diffusion process. This diffusion offers a novel interpretation of the Kalman-Bucy filter; i.e. as the conditional mean and covariance of an associated nonlinear McKean-Vlasov-type diffusion. This interpretation is interesting in its own right, and allows one to \emph{avoid Gaussian assumptions} on the initial state in a systematic way. We also introduce the concept of observability and controllability, which are signal related properties but which are relevant to the coming stability analysis. We also introduce the notion of a steady-state limit of the Riccati flow. 

 In Section \ref{sec-semigroup-props} we outline the relevant exponential and so-called Kalman-Bucy semigroups, associated largely with the Riccati flow and the stochastic flow of the Kalman-Bucy diffusion. We show how trajectories/solutions of the Riccati flow and Kalman-Bucy diffusion are defined in terms of these semigroups. Some preliminary technical lemmas are given concerning semigroup estimates and a number of invariance relationships are introduced. Both time-varying and homogeneous signal models are considered.

Section \ref{sec-riccati} is dedicated to the deterministic Riccati flow. Our first main result concerns the boundedness of the solution to the Riccati equation.
 
\begin{theo}\label{compact-localisation-intro}
There exists some $\upsilon>0$, assuming standard uniform observability and controllability conditions, and some 
$\Lambda_{min},\Lambda_{max}\in \SB^+_{r_1}$
such that
for any $t\geq \upsilon$  and  any $Q\in \SB_{r_1}^+$  we have
$$
\Lambda_{min}\leq \phi_{t}(Q)\leq \Lambda_{max}.
$$
\end{theo}

This result is stated precisely in Section \ref{sec-riccati} as Theorem \ref{compact-localisation} where the upper and lower bounds are given in terms of the observability and controllability Gramians. Indeed, this theorem correctly upper and lower bounds the Riccati flow in terms of these Gramians and it corrects Bucy's erroneous bounds given in \cite{bucy2}.

Section \ref{sec-riccati} is largely inspired by the seminal paper of Bucy \cite{bucy2} and we review, reprove, correct (where necessary), and refine those major results here. For example, under basic conditions, we consider the Lipschitz continuity and existence of solutions to the Riccati matrix differential equation. Following Bucy's original work \cite{bucy2}, a detailed proof of uniform convergence for the associated Kalman-Bucy semigroup is derived based on the corrected uniform bounds on the Riccati flow (and its inverse). This leads to a number of qualitative and \emph{quantitative} contraction estimates for the semigroup and the Riccati flow, both with time-varying models, and time-invariant models (where convergence to the fixed point of the Riccati operator is then considered). 

The first main result of this type is of the following form and is stated precisely in Section \ref{sec-riccati} as Theorem \ref{theo-expo-sg-ricc}.
 
\begin{theo}\label{theo-expo-sg-ricc-intro}
There exists some $\upsilon>0$ such that  for any $t\geq \upsilon$ and $Q_1,Q_2\in \SB^+_{r_1}$ we have
$$
\Vert\phi_t(Q_1)-\phi_t(Q_2)\Vert_2\leq 
\alpha\exp{\left\{-\beta t\right\}}~\Vert Q_1-Q_2\Vert_2
$$
for some positive parameters $(\alpha,\beta)$ whose values only depend on $(\Lambda_{min},\Lambda_{max})$.
The same inequality holds for any time $t\geq 0$ for some $\alpha=\alpha (Q_1,Q_2)$ that also depends on $(Q_1,Q_2)$.
\end{theo}

In Section \ref{sec-diffusion-flows} we initiate a novel analysis on the convergence of Kalman-Bucy stochastic flows, both in the classical filtering form, and the novel nonlinear diffusion form. The first main result of this type is the classical filtering stability result. 

\begin{theo} \label{conv-timevarying-KF-intro}
There exists some $\upsilon>0$ such that  for any $t\geq s\geq \upsilon$ it follows
$$
\sup_{Q\in\SB^+_{r_1}}\left\Vert\,  \EE\left( \psi_{s,t}(x,Q) - \varphi_{s,t}(X_s) \,\vert\, X_s\right)\, \right\Vert_2
\leq \alpha\exp{\left\{-\beta (t-s)\right\}}\Vert\,x-X_s \Vert_2
$$
with the parameters $\alpha,\beta >0$ as given in Theorem \ref{theo-expo-sg-ricc-intro}.

\end{theo}

This result is stated precisely in Section \ref{sec-diffusion-flows} as Theorem \ref{conv-timevarying-KF} and it is interesting because it shows that the bias between the filter and the signal is exponentially stable irregardless of the stability properties of the (time-varying) true signal. Much more is true, and we study exponential and comparison inequalities that bound with dedicated probability (at any time), the stochastic flow of the filter sample paths with respect to the underlying signal. That is, the next theorem shows that all the sample paths of the Kalman filter remain bounded close to the true signal with a large exponential probability. This result is stated precisely in Section \ref{sec-diffusion-flows} as Theorem \ref{prop-event-control-KB}.

\begin{theo}
The conditional probability of the following
events
$$
\begin{array}[t]{l}
\left\Vert~\psi_{s,t}(x,Q)-\varphi_{s,t}(X_s)\right\Vert_2 ~\leq~  \displaystyle\alpha_1(Q)~e^{-\beta(t-s)}~\Vert x-X_s\Vert_2+ \alpha_2(Q)\left[1+\delta+\sqrt{\delta}\right]
\end{array}
$$
given the state variable $X_s$ is greater than $1-e^{-\delta}$, for any $\delta\geq 0$ and any $t\in[s,\infty[$,
and some parameter $\beta >0$ and some $\alpha_i(Q)$ whose values only depend on $Q$, with $i=1,2$.
\end{theo}

In addition to this probabilistic convergence result, we give almost sure contraction-type estimates on the mean squared stochastic flow of the filter, conditioned on the underling signal of interest. The next result is stated precisely in Section \ref{sec-diffusion-flows} as Theorem \ref{theo-psi-st}.

\begin{theo}\label{theo-psi-st-intro}
For any $t\geq s\geq 0$, $x_1,x_2\in \RR^{r_1}$, $Q_1,Q_2\in\SB^+_{r_1}$ and $n\geq 1$ we have the almost sure local contraction estimate
 $$
\begin{array}{l}
\EE\left(\Vert\psi_{s,t}(x_1,Q_1)-\psi_{s,t}(x_2,Q_2)\Vert_2^{2n}~\vert~X_s\right)^{\frac{1}{2n}}
\\
\\
\qquad\leq \displaystyle \alpha_1(Q_1,Q_2)~e^{-\beta (t-s)}~\Vert x_1-x_2 \Vert_2+e^{-\beta(t-s)}~  \alpha_2(Q_1,Q_2)
\left\{\left\Vert x_2-X_s\right\Vert_2+\sqrt{n}\right\}~\Vert Q_1-Q_2\Vert_2
\end{array}$$
for some $\beta >0$ and some $\alpha_i(Q_1,Q_2)$ whose values only depend on $(Q_1,Q_2)$, with $i=1,2$.

\end{theo}

The preceding two results concern the Kalman-Bucy filter. We also have analogous results for the nonlinear Kalman-Bucy diffusion, i.e. we show that all sample paths of the Kalman-Bucy diffusion follow the true signal with a large exponential probability and we provide an almost sure contraction-type estimate on the mean squared stochastic flow of both diffusion.

\begin{theo}
The conditional probability of the following
events
$$
\begin{array}[t]{l}
\left\Vert~\overline{\psi}_{s,t}(x,Q)-\varphi_{s,t}(X_s)\right\Vert_2 ~\leq~  \displaystyle\alpha_1(Q)~e^{-\beta(t-s)}~\Vert x-X_s\Vert_2+ \alpha_2(Q)\left[1+\delta+\sqrt{\delta}\right]
\end{array}
$$
given the state variable $X_s$ is greater than $1-e^{-\delta}$, for any $\delta\geq 0$ and any $t\in[s,\infty[$,
and some parameter $\beta >0$ and some $\alpha_i(Q)$ whose values only depend on $Q$, with $i=1,2$.
\end{theo}

\begin{theo}
For any $t\geq s\geq 0$, $x_1,x_2\in \RR^{r_1}$, $Q_1,Q_2\in\SB^+_{r_1}$ and $n\geq 1$ we have
the almost sure local contraction estimate
 $$
\begin{array}{l}
\EE\left(\Vert \overline{\psi}_{s,t}(x_1,Q_1)-\overline{\psi}_{s,t}(x_2,Q_2)\Vert_2^{2n}~\vert~X_s\right)^{\frac{1}{2n}}
\\
\\
\qquad\leq \displaystyle \alpha_1(Q_1,Q_2)~e^{-\beta (t-s)}~\Vert x_1-x_2 \Vert_2+e^{-\beta(t-s)}~  \alpha_2(Q_1,Q_2)
\left\{\left\Vert x_2-X_s\right\Vert_2+\sqrt{n}\right\}~\Vert Q_1-Q_2\Vert_2
\end{array}$$
for some $\beta >0$ and some $\alpha_i(Q_1,Q_2)$ whose values only depend on $(Q_1,Q_2)$, with $i=1,2$.
\end{theo}

The preceding two results are stated precisely in Section \ref{sec-diffusion-flows} as Theorems \ref{prop-event-control-KB-diffusion} and \ref{theo-psi-diffusion-st} respectively. Both of the preceding results offer a general notion of filter stability. To the best of our knowledge, this approach to studying the stability of Kalman-Bucy stochastic flows is novel; and indeed this is certainly true so far as the nonlinear Kalman-Bucy diffusion is concerned.

Throughout this article, attention is paid to the quantitative nature of the convergence and stability results. For example, we study exponential rates in the autonomous case, in terms of different estimates on the relevant semigroups, and we track closely the related constants in front of these exponential terms. Our estimates are explicitly expressed with local Lipschitz contraction inequalities, dependent on the relevant signal matrix norms, etc. This contrasts with the classical analysis in \cite{kalman61,bucy68,anderson71,ocone-pardoux}, which is purely qualitative in nature. 

Carefully tracking constants is important for many applications; e.g. when studying the stability of ensemble Kalman filters \cite{legland09,tong16,dm-16-enkf}, extended Kalman filters \cite{reif2000,dm-16-eekf}, or when it comes to understanding general approximations of the Kalman filter, and it's error dependence on the state-space dimension \cite{mf-dm-04,rebeschini2015,tong-16-enkf}.

\subsection{Some basic notation}
This section details some basic notation and terms used throughout the article.

 With a slight abuse of notation, we denote by $Id$ the identity matrix (with the size obvious from the context). The matrix transpose is denoted by $'$.

Denote by $\lambda_i(A)$, with $1\leq i\leq r$, the non-increasing sequence of eigenvalues of a $(r\times r)$-matrix $A$ and let $\mbox{\rm Spec}(A)$ be the set of all eigenvalues. We denote by $\lambda_{min}(A)=\lambda_{r}(A)$ and $\lambda_{max}(A)=\lambda_{1}(A)$ the minimal and the maximal eigenvalue. Set $A_{sym}:=(A+A^{\prime})/2$ for any $(r\times r)$-square matrix $A$. 

We define a particular logarithmic norm $\mu(A)$ of an $(r_1\times r_1)$-square matrix $A$ by
\begin{equation}\label{def-log-norm}
\begin{array}{rcl}
\mu(A)&:=&\inf{\{\alpha:\forall x,~\langle x,Ax\rangle\leq \alpha \left\Vert x\right\Vert_2^2\}} \\
&=&\lambda_{max}\left(A_{sym}\right)\\
&=&\inf{\{\alpha:\forall t\geq 0,~\Vert \exp{(At)}\Vert_2\leq \exp{(\alpha t)}\}}.
\end{array}
\end{equation}
The above equivalent formulations show that 
$$
\mu(A)~\geq~ \varsigma(A):=\max{\left\{\mbox{\rm Re}(\lambda)~:~\lambda\in \mbox{\rm Spec}(A)\right\}}
$$
where $\mbox{\rm Re}(\lambda)$ stands for the real part of the eigenvalues $\lambda$.
The parameter $\varsigma(A)$ is often called the spectral abscissa of $A$. Also note that $A_{sym}$ is negative definite as soon as $\mu(A)<0$. 

The Frobenius matrix norm of a given $(r_1\times r_2)$ matrix $A$ is defined by
$$
\left\Vert A\right\Vert_{F}^2=\mbox{\rm tr}(A^{\prime}A)
\qquad\mbox{\rm
with the trace operator $\mbox{\rm tr}(\cdot)$.}
$$
If $A$ is a matrix $(r\times r)$, we have $\left\Vert A\right\Vert_{F}^2=\sum_{1\leq i,j\leq r}A(i,j)^2\geq \Vert A\Vert^2$. For any $(r\times r)$-matrix $A$, we recall norm equivalence formulae
$$\Vert A\Vert_2^2=\lambda_{max}(A^{\prime}A)\leq \mbox{\rm tr}(A^{\prime}A)=\Vert A\Vert_F^2\leq r~\Vert A\Vert_2^2$$.

The Hoffmann-Wieland theorem (Theorem 9.21 in~\cite{fielder}) also tells us that for any symmetric matrices $A,B$ we have
$$
\sum_{1\leq i\leq r_1}\left(\lambda_i(A)-\lambda_i(B)\right)^2\leq \Vert A-B\Vert_F^2=\sum_{1\leq i\leq r_1}\left(\lambda_i(A-B)\right)^2.
$$

Now, given some random variable $Z$ with some probability measure or distribution $\eta$ and some measurable function
$f$ on some product space $\RR^r$, we let $$\eta(f)=\EE(f(Z))=\int~f(x)~\eta(dx)$$ be the integral of $f$ w.r.t. $\eta$ or the expectation of $f(Z)$.  As a rule any multivariate variable, say $Z$,
is represented by a column vector and we use the transposition operator $Z^{\prime}$ to denote the row vector (similarly for matrices).

We also recall as background that for any non-negative random variable $Z$ such that
$$
\EE\left(Z^{2n}\right)^{1/n}\leq z^2~n\quad\mbox{\rm for some parameter $z\not=0$}
$$
and for any $n\geq 1$ we have
$$
\EE\left(Z^{2n}\right)\leq (z^2n)^n\leq \frac{e}{\sqrt{2}}~\left(\frac{e}{2}~z^2\right)^{n} \EE(V^{2n})
$$
for some Gaussian and centered random variable $V$ with unit variance. We check this claim using Stirling approximation
\begin{eqnarray*}
\EE(V^{2n})&=&2^{-n}\frac{(2n)!}{n!}\nonumber\\
&\geq& e^{-1}~2^{-n}\frac{\sqrt{4\pi n}~(2n)^{2n}~e^{-2n}}{\sqrt{2\pi n}~n^{n}~e^{-n}}=\sqrt{2}e^{-1}~\left(\frac{2}{e}\right)^{n}~~n^n. 
\end{eqnarray*}
By \cite[Proposition 11.6.6]{mf-dm-13}, the probability of the following event
\begin{equation}\label{event-control}
(Z/z)^2\leq \frac{e^2}{\sqrt{2}}~\left[\frac{1}{2}+\left(\delta+\sqrt{\delta}\right)\right]
\end{equation}
is greater than $1-e^{-\delta}$, for any $\delta\geq 0$.

Given a real valued continuous martingale $M_t$ starting at the origin $M_0=0$, for any $n\geq 1$ and any 
time horizon $t\geq 0$ we have
\begin{equation}\label{BDG}
\EE\left(\sup_{0\leq s\leq t}\vert M_s\vert^n\right)^{1/n}\leq 2\sqrt{2}~\sqrt{n}~\EE\left(\langle M\rangle_t^{n/2}\right)^{1/n}.
\end{equation}
Proof of these Burkholder-Davis-Gundy inequalities is in~\cite{yao-feng-reng}; see also~\cite[Theorem B.1 p. 97]{Khoshnevisan}.

\section{Description of the models}\label{desc-sec-intro}

\subsection{The Kalman-Bucy filter}

In general (i.e. not assuming $X_0$ is Gaussian), for any $0\leq s\leq t$, we define the stochastic flow 
$$
  \Phi_{s,t}~:~
 (x,Q)\in (\RR^{r_1}\times \SB^+_{r_1})\mapsto \Phi_{s,t}(x,Q)=\left(\psi_{s,t}(x,Q),\phi_{s,t}(Q)\right)\in (\RR^{r_1}\times \SB^+_{r_1})
 $$
as describing the Kalman-Bucy filter; where for any horizon $s$ and any time $t\in [s,\infty[$ we have
\begin{equation*}
\left\{
\begin{array}{rcl}
d\psi_{s,t}(x,Q)&=&\left[A_t-\phi_{s,t}(Q)S_t\right]~\psi_{s,t}(x,Q)~dt+\phi_{s,t}(Q)~C_t^{\prime}R^{-1}_{2}~dY_t\\
&&\\
\partial_t\phi_{s,t}(Q)&=&\ricc\left(\phi_{s,t}(Q)\right)
\quad\mbox{\rm with}\quad\Phi_{s,s}(x,Q)=(x,Q).
\end{array}
\right.
\end{equation*}
With similar notation, we also denote by $\varphi_{s,t}(x)$ the stochastic flow of the signal process,
$$
d\varphi_{s,t}(x)=A_t\varphi_{s,t}(x)dt+R_1^{1/2}~dW_t\quad\mbox{\rm with}\quad
\varphi_{s,s}(x)=x
$$
for any $t\in [s,\infty[$, and any $x\in \RR^{r_1}$.

Note that in general $\phi_{s,t}(\phi_s(Q))=\phi_t(Q)$ and, 
$$
	\phi_{s+t}(Q) = \phi_{s,s+t}(\phi_{s}(Q))~~\mathrm{and}~~\phi_{s,t}(Q) = \phi_{s+u,t}(\phi_{s,s+u}(Q)), ~ 0\leq u\leq t-s.
$$
Observe that when the signal is time-invariant, then so is the Riccati equation and thus 
$$
	\phi_{s,s+t}(Q)=\phi_{t}(Q)=:\phi_{0,t}(Q) ~~\mathrm{or}~~\phi_{s,t}(Q) = \phi_{t-s}(Q)=\phi_{u,t-s}(\phi_{u}(Q)), ~ 0\leq u\leq t-s
$$
along with numerous other (equivalent) combinations of stationary shifts.

\subsection{Nonlinear Kalman-Bucy diffusions}

For any $0\leq s\leq t$, we let 
$$
  \overline{\Phi}_{s,t}~:~
 (x,Q)\in (\RR^{r_1}\times \SB^+_{r_1})\mapsto  \overline{\Phi}_{s,t}(x,Q)=\left(\overline{\psi}_{s,t}(x,Q),\phi_{s,t}(Q)\right)\in (\RR^{r_1}\times \SB^+_{r_1})
 $$
 be the stochastic flow of the Kalman-Bucy diffusion; that is, for any time horizon
$s$ and any time $t\in [s,\infty[$ we have
\begin{eqnarray*}
d\overline{\psi}_{s,t}(x,Q)&=&\left[A_t-\phi_{s,t}(Q)S_t\right]~\overline{\psi}_{s,t}(x,Q)~dt+\phi_{s,t}(Q)~C_t^{\prime}R^{-1}_{2}~dY_t\\
&&\hskip5cm+R^{1/2}_{1}~d\overline{W}_t-\phi_{s,t}(Q)C_t^{\prime}R^{-1/2}_{2}d\overline{V}_{t}
\end{eqnarray*}
with $\overline{\psi}_{s,s}(x,Q)=x$, for $t=s$.

\subsection{Observability and controllability conditions}

We consider the observability and controllability Gramians $\Oa_{s,t}$ and $\Ca_{s,t}$ defined by
$$
 \Ca_{s,t}:=\int_{s}^{t}~\exp{\left[\oint_r^{t} A_u~du\right]}R_1\exp{\left[\oint_r^{t} A_u~du\right]}'~dr
$$
and
$$
 \Oa_{s,t}:= \int_{s}^{t}~\exp{\left[\oint_{t}^{r} A_u~du\right]}'S_r\exp{\left[\oint_{t}^{r} A_u~du\right]}~dr
$$
for all $t\geq s\geq0$. We let $\Ca_{t}:=\Ca_{0,t}$ and $\Oa_{t}:=\Oa_{0,t}$. Here $\exp{[\oint_s^{t} A_u~du]}$ defines a semigroup associated with the matrix flow. We return to this semigroup and its properties later; see (\ref{eq-semigroup-def-good}) for a more precise definition. 

We take the following assumption as holding in the statement of all results. 

\textbf{Standing Assumption:} {\emph{
There exists parameters $\upsilon,\varpi^{o,c}_{\pm}>0$ such that
\begin{equation}\label{steady-state-eq-2}
\varpi_-^c~Id\leq \Ca_{t,t+\upsilon}
\leq \varpi_+^c~Id
\quad\mbox{\rm
and}
\quad
\varpi_-^o~Id\leq \Oa_{t,t+\upsilon}\leq \varpi_+^o~Id
\end{equation}
uniformly for all $t\geq0$. The parameter $\upsilon$ is called the 
interval of observability/controllability. }}

 Note that if the signal matrices are time-invariant, then the pair $(A,R_1^{1/2})$ is a controllable and $(A,C)$ is observable if \begin{equation}\label{def-contr-obs}
\left[R_1^{1/2},~AR_1^{1/2},~\ldots, A^{r_1-1}R_1^{1/2}\right]\quad
\mbox{\rm and}\quad
\left[\begin{array}{c}
C\\
CA\\
\vdots\\
CA^{r_1-1}
\end{array}
\right]
 \end{equation}
both have rank $r_1$. Under our common (in filtering) modelling assumption that $R_1\in\SB_{r_1}^+$, the pair $(A,R_1^{1/2})$ is immediately controllable. 

Under these conditions (\ref{def-contr-obs}) there always exists parameters $\upsilon,\varpi^{o,c}_{\pm}>0$ ensuring that (\ref{steady-state-eq-2}) holds. For example, whenever the signal drift matrix $A$ is diagonalizable, and $R_1,S\in\SB_{r_1}^+$ we can choose
$$
\varpi_-^c=\lambda_{min}(R_1)~\min_{\lambda\in\mbox{\rm Spec}(A)}\frac{e^{2\lambda\upsilon}-1}{2\lambda}~\leq \varpi_+^c=
\lambda_{max}(R_1)~\max_{\lambda\in\mbox{\rm Spec}(A)}\frac{e^{2\lambda\upsilon}-1}{2\lambda}~
$$
as well as
$$
\varpi_-^o=\lambda_{min}(S)~\min_{\lambda\in\mbox{\rm Spec}(A)}\frac{1-e^{-2\lambda\upsilon}}{2\lambda}~\leq \varpi_+^o=
\lambda_{max}(S)~\max_{\lambda\in\mbox{\rm Spec}(A)}\frac{1-e^{-2\lambda\upsilon}}{2\lambda}~
$$
for any $\upsilon>0$.

In the time-invariant case, these conditions (\ref{def-contr-obs}) are sufficient (but not necessary \cite{Poubelle86,callier95}) to ensure there exists a (unique) positive definite fixed-point matrix $P=\phi_t(P)$ solving the so-called algebraic Riccati equation
\begin{equation}\label{steady-state-eq}
\ricc(P):=AP+PA^{\prime}-PSP+R_1=0.
\end{equation}
In this time-invariant model setting, the matrix difference $A-PS$ is asymptotically stable (Hurwitz stable) even when the signal matrix $A$ is unstable \cite[Theorems 9.12, 9.15]{Lancaster1995}. More relaxed conditions (i.e. detectability and stabilisability) for a stabilising solution (perhaps only positive semi-definite) to exist are discussed widely in the literature; see \cite{kucera72,Molinari77,Lancaster1995} and the convergence results in \cite{Kwakernaak72,callier81}. A (marginally) stable solution to (\ref{steady-state-eq}) exists under a detectability condition, and convergence to this solution is given under mild additional conditions in \cite{Poubelle86,callier95,Park97}. We also note that the stability of $A-PS$ follows from Theorems \ref{theo-expo-sg-ricc} and \ref{conv-timevarying-KF}; and it follows that $\phi_{t}(Q)$ for $Q\in \SB_{r_1}^+$ converges to the fixed point $P$ due to Theorem \ref{theo-expo-sg-ricc} and the Banach fixed-point theorem.

In the time-varying case, we are interested in asymptotic stability results and results that tend to bound the Riccati flow $\phi_{t}(Q)$ uniformly on both sides by the controllability and observability Gramians. In this setting, there is typically no fixed point for the flow $\phi_{t}(Q)$ and the difference $A_t-\phi_{t}(Q)S_t$ need not be a stable matrix at any instant in general; see also \cite{gevers85,Poubelle88}.

Since we switch between time-invariant signal models and time-varying models (in which no fixed point of (\ref{nonlinear-KB-Riccati}) generally exists), we choose not to relax our observability/controllability assumptions, e.g. viz. \cite{kucera72,Poubelle86}. In fact, as discussed in \cite{callier95}, there is really no generality lost by assuming observability over the (arguably) weaker detectability condition \cite{Poubelle86} even in the time-invariant case. One may substitute a generalised controllability condition like that studied in \cite{anderson71}; albeit only weaker results are achievable (as shown by example in \cite{anderson71}). Anyway, in our case, we have $R_1\in\SB_{r_1}^+$, as is common in many filtering problems.

\section{Semigroups of the Riccati flow and Kalman-Bucy filter}\label{sec-semigroup-props}

\subsection{Exponential semigroups}

The transition matrix associated with a smooth flow of $(r\times r)$-matrices $A:u\mapsto A_u$ is denoted by
\begin{equation}\label{eq-semigroup-def-good}
\Ea_{s,t}(A)=\exp{\left[\oint_s^t A_u~du\right]}\Longleftrightarrow \partial_t \Ea_{s,t}(A)=A_t~\Ea_{s,t}(A)\quad\mbox{\rm and}\quad
\partial_s \Ea_{s,t}(A)=-\Ea_{s,t}(A)~A_s
\end{equation}
for any $s\leq t$, with $\Ea_{s,s}=Id$, the identity matrix. 

The following technical lemma gives a pair of semigroup estimates for the state transition matrices associated with a sum of drift-type matrices. 
\begin{lem}\label{perturbation-lemma-intro}
Let $A:u\mapsto A_u$ and $B:u\mapsto B_u$ be the smooth flows of $(r\times r)$-matrices.
For any $s\leq t$ we have
  $$
 \left\Vert   \Ea_{s,t}(A+B)\right\Vert_2\leq \exp{\left(\int_s^t\mu(A_u)~du+\int_s^t~\Vert B_u\Vert_2~du\right)}.
 $$
In addition, for the matrix spectral, or Frobenius, norm $\Vert\cdot\Vert$ we have
 $$
  \left\Vert   \Ea_{s,t}(A+B)\right\Vert\leq \alpha_A \exp{\left[-\beta_A(t-s)+\alpha_A\int_s^t \Vert B_u\Vert~du\right]}
 $$
 as soon as
  $$
 \forall 0\leq s\leq t\qquad \Vert  \Ea_{s,t}(A) \Vert\leq \alpha_A~\exp{\left(-\beta_A~(t-s)\right)}.
 $$

\end{lem}
\proof
The above estimate is a direct consequence of the matrix log-norm inequality
$$
\mu(A_t+B_t) \leq \mu(A_t)+\mu(B_t)\quad\mbox{\rm and the fact that}\quad
\mu(B_t)\leq \Vert B_t\Vert_2.
$$
This ends the proof of the first assertion. To check the second assertion we observe that
$$
\partial_{t} \Ea_{s,t}(A+B) =\left(\partial_{t} \Ea_{t}(A+B)\right) \Ea_{s}(A+B)^{-1}=A_t\Ea_{s,t}(A+B)+B_t\Ea_{s,t}(A+B).
$$
This implies that
 $$
 \Ea_{s,t}(A+B) =
\Ea_{s,t}(A)+\int_s^t~\Ea_{u,t}(A)~B_u\Ea_{s,u}(A+B)~du
 $$
 for any $s\leq t$ from which we prove that
   \begin{eqnarray*}
 e^{\beta_A(t-s)}\Vert\Ea_{s,t}(A+B) \Vert&\leq &
\alpha_A+\alpha_A~\int_s^t~ e^{\beta_A(t-s)}~ e^{-\beta_A(t-u)}~\Vert B_u\Vert~\Vert
\Ea_{s,u}(A+B)\Vert~du\\
&=&\alpha_A+\alpha_A~\int_s^t~ 
\Vert B_u\Vert~ e^{\beta_A(u-s)}\Vert
\Ea_{s,u}(A+B)\Vert~du
\end{eqnarray*}
By Gr\"onwall's lemma this implies that
$$
  e^{\beta_A(t-s)}\Vert\Ea_{s,t}(A+B) \Vert\leq \alpha_A \exp{\left[\int_s^t \alpha_A\Vert B_u\Vert~du\right]}.
$$
This ends the proof of the lemma.
\qed

\subsubsection{Time-invariant exponential semigroups}

For time-invariant matrices $A_t=A$, the state transition matrix reduces to a conventional matrix exponential $$\Ea_{s,t}(A)=e^{(t-s)A}= \Ea_{t-s}(A).$$ 

\vspace{0.1cm}In this subsection we are interested in estimating the norm of $\Ea_t(A)$. We state the following general convergence result on the time-invariant semigroup generated by the matrix $(A-PS)$ where $P$ is the fixed-point solution to (\ref{steady-state-eq}).

\begin{lem}\label{cor-hypothesis}
	 Under the time-invariant observability/controllability conditions (\ref{def-contr-obs}), it follows that,
\begin{equation}\label{hypothesis}
	\exists \nu>0,~ \exists \kappa<\infty~:~\forall t\geq 0,\quad \Vert e^{t(A-PS)}\Vert_2~\leq~ \kappa\,e^{-\nu\, t}.
\end{equation}
\end{lem}

\proof
The observability/controllability rank conditions (\ref{def-contr-obs}) are sufficient to ensure the existence of a (unique) positive definite solution $P$ of (\ref{steady-state-eq}) and that $\varsigma(A-PS)<0$; see \cite{kucera72,Lancaster1995}. We know
$$
	\Vert e^{t(A-PS)}\Vert_2 ~\leq~ e^{\mu(A-PS)t}
$$
applies whenever $\mu(A-PS)<0$. Otherwise, we can also use any of the estimates presented below in (\ref{first-estimate-Jordan}), (\ref{first-estimate-Jordan-diag}), (\ref{estimate-Schur}), (\ref{common-estimate-Jordan-Schur}). 
\qed

The norm of $\Ea_t(A)$ can be estimated in various ways: The first is based on the Jordan decomposition
$T^{-1}AT=J$ of the matrix $A$ in terms of $k$ Jordan blocks associated with the eigenvalues
with multiplicities $m_i$, with $1\leq i\leq k$. In this situation, we have the Jordan type estimate
\begin{equation}\label{first-estimate-Jordan}
	e^{\varsigma(A)t}~\leq~ \Vert \Ea_t(A)\Vert_2\leq \kappa_{{Jor},t}(T) \,e^{\varsigma(A)t}
\end{equation}
with
$$
	\kappa_{{Jor},t}(T):=\left(\max_{0\leq j<n}~\frac{t^j}{j!}\right)~ \Vert T\Vert_2\Vert T^{-1}\Vert_2 \quad \mbox{\rm and}\quad n:=\max_{1\leq i\leq k}~ m_i.
$$
Note that $\kappa_{{Jor},t}(T)$ depends on the time horizon $t$ as soon as $A$ is not of full rank. In addition, whenever $A$ is close to singular, the condition number $\mbox{\rm cond}(T):=\Vert T\Vert_2\Vert T^{-1}\Vert_2$ tends to be very large. When $A$ is diagonalizable the above estimate becomes
\begin{equation}\label{first-estimate-Jordan-diag}
	e^{\varsigma(A)t}~\leq~ \Vert \Ea_t(A)\Vert_2 ~\leq~ \mbox{\rm cond}(T) \,e^{\varsigma(A)t}.
\end{equation}

Another method is based on the Schur decomposition $U^{\prime}AU=D+T$ in terms of an unitary matrix $U$, with
$D=\mbox{\rm diag}(\lambda_1(A),\ldots;\lambda_r(A))$ and a strictly triangular matrix $T$ s.t. $T_{i,j}=0$ for any $i\geq j$. In this case we have the Schur type estimate
\begin{equation}\label{estimate-Schur}
	\Vert \Ea_t(A)\Vert_2~\leq~  \kappa_{{Sch},t}(T)\,e^{\varsigma(A)t}
		~\quad \mbox{\rm with}\quad \kappa_{{Sch},t}(T):=\sum_{0\leq i\leq r}\frac{(\Vert T\Vert_2 t)^i}{i!}.
\end{equation}
The proof of these estimates can be found in \cite{vanloan,vanloan-19}. In both cases for any $\epsilon\in ]0,1]$ and  any $t\geq 0$ we have
\begin{equation}\label{common-estimate-Jordan-Schur}
	e^{\varsigma(A)t}~\leq~\Vert\Ea_t(A)\Vert_2\leq \kappa_A(\epsilon)\,e^{(1-\epsilon)\varsigma(A)t}~
\end{equation}
for some constants $\kappa_A(\epsilon)$ whose values only depend on the parameters $\epsilon$. When $A$
is  asymptotically stable; that is all its eigenvalues have negative real parts, for any positive definite matrix $B$ we have
$$
	e^{\varsigma(A)t}~\leq~\Vert\Ea_t(A)\Vert_2 ~\leq~ \mbox{\rm cond}(T)\,\exp{\left[-t/\Vert B^{-1/2} T~B^{-1/2}\Vert_2\right]}
$$
 with the positive definite matrix
$$
	T=\int_0^{\infty}~e^{A^{\prime}t}~B~e^{At}~dt\Longleftrightarrow A^{\prime}T+TA=-B.
$$
The proof of these estimates can be found in \cite[see e.g. Theorem 13.6]{kresemir}.

\subsection{Kalman-Bucy semigroups}

For any $s\leq t$ and $Q_1,Q_2\in\SB_{r_1}^+$ we set
$$
	E_{s,t}(Q_1,Q_2):=\exp{\left[\oint_s^t\left(A_u-\frac{\phi_{u}(Q_1)+\phi_{u}(Q_2)}{2}~S_u\right)~du\right]}
		\quad\mbox{and}\quad E_{s,t}(Q_1):=E_{s,t}(Q_1,Q_1).
$$
When $s=0$ we write $E_{t}(Q_1)$ and $E_{t}(Q_1,Q_2)$ in place of $E_{0,t}(Q_1)$ and $E_{0,t}(Q_1,Q_2)$. In this notation we have
$$
E_{s,t}(Q_1,Q_2)=E_{t}(Q_1,Q_2)E_{s}(Q_1,Q_2)^{-1}\quad\mbox{\rm and}\quad
E_{s,t}(Q_1)=E_{t}(Q_1)E_{s}(Q_1)^{-1}.
$$
We have the following important result.

\begin{prop}\label{proposition-phi-sg}
For any $s\leq t$ and $Q_1,Q_2\in\SB_{r_1}^+$ we have
\begin{eqnarray}
\phi_t(Q_1)-\phi_t(Q_2)
\displaystyle&=&~~E_{s,t}(Q_1)~~~~\left[\phi_s(Q_1)-\phi_s(Q_2)\right]~~~E_{s,t}(Q_2)^{\prime}\label{polarization-formulae-phi1}\\
&=&\displaystyle E_{s,t}(Q_1,Q_2)~\left[\phi_s(Q_1)-\phi_s(Q_2)\right]~E_{s,t}(Q_1,Q_2)^{\prime}\label{polarization-formulae-phi2}
\end{eqnarray}
as well as 
\begin{eqnarray}\label{polarization-formulae-phi3}
\phi_t(Q_1)-\phi_t(Q_2)&
\displaystyle=&E_{s,t}(Q_2)\left[\phi_s(Q_1)-\phi_s(Q_2)\right]~E_{s,t}(Q_2)^{\prime}\\
&&\hskip-.5cm\displaystyle-\int_s^t~
E_{u,t}(Q_2)~\left[\phi_u(Q_1)-\phi_u(Q_2)\right]~S_u~\left[\phi_u(Q_1)-\phi_u(Q_2)\right]~
E_{u,t}(Q_2)^{\prime}~du. \nonumber
\end{eqnarray}
\end{prop}

\proof
These semigroup formulae are direct consequences of the following 
three polarization-type formulae 
\begin{equation}\label{polarization-formulae}
\begin{array}{l}
\ricc(Q_1)-\ricc(Q_2) \qquad\\
\\
\qquad\qquad=(A_t-Q_1S_t)(Q_1-Q_2)+(Q_1-Q_2)(A_t-Q_2S_t)^{\prime}\\
\\
\qquad\qquad=\left[A_t-\frac{1}{2}(Q_1+Q_2)S_t\right](Q_1-Q_2)+(Q_1-Q_2)\left[A_t-\frac{1}{2}(Q_1+Q_2)S_t\right]^{\prime}\\
\\
\qquad\qquad=(A_t-Q_2S_t)(Q_1-Q_2)+(Q_1-Q_2)(A_t-Q_2S_t)^{\prime}-(Q_1-Q_2)S_t(Q_1-Q_2)\\
\end{array}
\end{equation}
where the first line implies (\ref{polarization-formulae-phi1}), the second line implies (\ref{polarization-formulae-phi2}), and the third line implies (\ref{polarization-formulae-phi3}). We check these polarization-type formulae using the decompositions
\begin{eqnarray*}
Q_1S_tQ_1-Q_2S_tQ_2&=&Q_1S_t(Q_1-Q_2)+(Q_1-Q_2)S_tQ_2\\
&=&\frac{1}{2}(Q_1+Q_2)S_t(Q_1-Q_2)+\frac{1}{2}(Q_1-Q_2)S_t(Q_1+Q_2)\\
&=&(Q_1-Q_2)S_t(Q_1-Q_2)+Q_2S(Q_1-Q_2)+(Q_1-Q_2)S_tQ_2.
\end{eqnarray*}
The proof of (\ref{polarization-formulae-phi1}), (\ref{polarization-formulae-phi2}) and (\ref{polarization-formulae-phi3}) from basic calculations, for example, equation (\ref{polarization-formulae-phi2}) follows by
\begin{eqnarray}
\partial_t(\phi_t(Q_1)-\phi_t(Q_2))&=& \ricc(\phi_t(Q_1))-\ricc(\phi_t(Q_2))\nonumber \\
&= & \left(A_t-\frac{\phi_t(Q_1)+\phi_t(Q_2)}{2}S_t\right)(\phi_t(Q_1)-\phi_t(Q_2)) \nonumber \\
&& \qquad\qquad\qquad + ~(\phi_t(Q_1)-\phi_t(Q_2))\left(A_t-\frac{\phi_t(Q_1)+\phi_t(Q_2)}{2}S_t\right)^{\prime}. \nonumber 
\end{eqnarray}
Now the solution of this linear equation is given by $E_{s,t}(Q_1,Q_2)~\left[\phi_s(Q_1)-\phi_s(Q_2)\right]~E_{s,t}(Q_1,Q_2)^{\prime}$ which is (\ref{polarization-formulae-phi2}). This ends the proof of the proposition. \qed

Given a time-varying signal model, it is useful to define some additional notation. For any $s\leq u\leq t$ and $Q\in\SB_{r_1}^+$ we set
$$
E_{t\vert s}(Q) := \exp{\left[\oint_s^t\left(A_u-\phi_{s,u}(Q)S_u\right)~du\right]}
	\quad\mbox{and}\quad
		E_{u,t\vert s}(Q):= \exp{\left[\oint_u^t\left(A_r-\phi_{s,r}(Q)S_r\right)~dr\right]} 
$$
with $E_{u,t\vert s}(Q)=E_{t\vert s}(Q)E_{u\vert s}(Q)^{-1}$. Note there is a relationship between $E_{t\vert s}$ and $E_{s,t}$ in the following sense
\begin{eqnarray*}
	E_{t\vert s}(\phi_s(Q)) &=& \exp{\left[\oint_s^{t}\left(A_u-\phi_{s,u}(\phi_s(Q))S_u\right)du\right]}\\
		&=& \exp{\left[\oint_s^{t}\left(A_u-\phi_{u}(Q)S_u\right)du\right]} ~=~ E_{s,t}(Q)
\end{eqnarray*}
where we simply used $\phi_{s,t}(\phi_s(Q))=\phi_{t}(Q)$.

In the time-invariant signal model setting, we point to \cite{bd-CARE} for an explicit expression of this exponential semigroup $E_{s,t}(Q)$, along with some applications in a refined stability analysis of the associated Riccati equation. We return to this latter remark later.

\section{Riccati flows} \label{sec-riccati}

We start this section with a preliminary result concerning the monotonicity of the Riccati operator, some basic boundedness results and a Lipschitz estimate.

\begin{prop}
The Riccati flow $Q\mapsto \phi_t(Q)$ is a non-decreasing function w.r.t. the Loewner partial order; that is we have
$$
Q_1\leq Q_2\Longleftrightarrow\phi_t(Q_1)\leq\phi_t(Q_2).
$$
For any $Q_1,Q_2\in\SB^+_{r_1}$ we have the local Lipschitz inequality
\begin{equation}\label{lipschitz-phi}
\Vert \phi_t(Q_1)-\phi_t(Q_2)\Vert_F\leq l_{Q_1,Q_2}( \phi_t)~\Vert Q_1-Q_2\Vert_F
\end{equation}
for some Lipschitz constant 
$$
l_{Q_1,Q_2}( \phi_t)\leq \left[\left\Vert E_t(Q_1)\right\Vert_2\,\left\Vert E_t(Q_2)\right\Vert_2\right]\wedge \left\Vert E_t(Q_1,Q_2)\right\Vert_2<\infty.
$$
\end{prop}
\proof
Using Proposition~\ref{proposition-phi-sg} we prove that $Q\mapsto \phi_t(Q)$ is an non decreasing function w.r.t. the Loewner partial order. The Lipschitz estimate (\ref{lipschitz-phi}) is a direct consequence of 
the implicit semigroup formulae (\ref{polarization-formulae-phi1}) and (\ref{polarization-formulae-phi2}).
\qed

It follows that for any $Q\in\SB^+_{r_1}$, the time-varying Riccati flow $\phi_t(Q)$ is well-defined and a unique solution exists for all $t\geq0$; since the local Lipschitz estimate is `global' on any finite interval. 

The Riccati flow $Q\mapsto \phi_t(Q)$ also depends monotonically on the parameters $S$ and $R_1$. 

\begin{cor}
Let $R_1(2)\geq R_1(1)\in \SB_{r_1}^+$ and $S_t(1)\geq S_t(2)\geq0$ for all $t\geq0$. Then $\phi_{s,t}(Q,2)\geq\phi_{s,t}(Q,1)$ where 
$$\partial_t\phi_{s,t}(Q,i)=A_tQ+QA_t^{\prime}-QS_t(i)Q+R_1(i),~~i\in\{1,2\}.$$
\end{cor}

Now define
\begin{equation}\label{ref-phi-uniform-max}
	\Vert\phi(Q)\Vert_2:=\sup_{t\geq 0}\Vert \phi_t(Q)\Vert_2 <\infty
\end{equation}
which is always uniformly bounded for small enough $t$ as a result of the Lipschitz estimate on $[0,t]$. 

In the time-invariant setting, when the desired solution $P$ of (\ref{steady-state-eq}) exists, the following result characterizes a uniform upper-bound on the Riccati flow and a bound on its growth.

\begin{prop}\label{prop-lipschitz-continuity}
The Riccati flow obeys
\begin{equation*}
	Q\mapsto \phi_t(Q)\leq P+E_t(P) (Q-P)E_t(P)^{\prime}.
\end{equation*} 
In addition, for any $Q\in\SB^+_{r_1}$ we have the uniform estimates
\begin{equation}\label{ref-phi-maxQ}
	\Vert\phi(Q)\Vert_2 \leq \Vert P\Vert_2+\kappa^2\Vert Q-P\Vert_2~\quad\mathrm{and}\quad~\sup_{t>0}~t^{-1}\Vert \phi_t(Q)-Q\Vert_2<\infty
\end{equation} 
where $\kappa$ is defined in Lemma \ref{cor-hypothesis}.
\end{prop}

\proof
Choosing $Q_2=P$ and $s=0$ in (\ref{polarization-formulae-phi3})
we find that
\begin{eqnarray*}
	\phi_t(Q)-P=E_{t}(P)\left[Q-P\right]E_{t}(P)^{\prime}-\int_0^t\,E_{u,t}(P)\,\left[\phi_u(Q)-P\right]\,S\,\left[\phi_u(Q)-P\right]\,E_{u,t}(P)^{\prime}~du.
\end{eqnarray*}
This implies
\begin{eqnarray*}
	0\leq \phi_t(Q)\leq P+E_t(P) (Q-P)E_t(P)^{\prime}&\Rightarrow &\Vert \phi_t(Q)\Vert_2\leq \Vert P\Vert_2+\Vert E_t(P) \Vert^2_2\Vert Q-P\Vert_2.
\end{eqnarray*} 
It then follows that $\Vert\phi(Q)\Vert_2 \leq \Vert P\Vert_2+\kappa^2\Vert Q-P\Vert_2$ from which we conclude
\begin{eqnarray*}
	\phi_t(Q)=Q+\int_0^t~\mbox{\rm Ricc}(\phi_s(Q))~ds\Rightarrow \Vert \phi_t(Q)-Q\Vert_F\leq c_Q~t
\end{eqnarray*}
for some finite constant $c_Q$ whose values only depends on $Q$. This completes the proof.
\qed

\subsection{Uniform bounds on the Riccati flow}

We let $\Ca_t(\Oa)$ and $\Oa_t(\Ca)$ be the Gramian matrices defined by
\begin{eqnarray*}
\Oa_t(\Ca)&:=&\Ca_t^{-1}\left[\int_0^t~\Ea_{s,t}(A)~\Ca_s~ S_s~\Ca_s~\Ea_{s,t}(A)^\prime~ds\right]\Ca_t^{-1}\\
\Ca_t(\Oa)&:=& \Oa_t^{-1}\left[\int_0^t~\Ea_{s,t}^\prime(A)^{-1}~\Oa_sR_1~\Oa_s~\Ea_{s,t}(A)^{-1}~ds\right]\Oa_t^{-1}.
\end{eqnarray*}
Under our standard observability and controllability assumptions (\ref{steady-state-eq-2}), there exists some parameters $\varpi^{c}_{\pm}(\Oa),\varpi^{o}_{\pm}(\Ca)>0$ such that
$$
 \varpi_-^c(\Oa)~Id
\leq \Ca_{\upsilon}(\Oa)\leq \varpi_+^c(\Oa)~Id\qquad
\mbox{\rm and}\qquad
 \varpi_-^o(\Ca)~Id
\leq \Oa_{\upsilon}(\Ca)\leq \varpi_+^o(\Ca)~Id
$$
holds uniformly on the interval $\upsilon>0$ of observability/controllability.

The main objective of this section is to prove the following theorem.

\begin{theo}\label{compact-localisation}
For any $t\geq \upsilon$  and  any $Q\in \SB_{r_1}^+$  we have
$$
\left(\Oa_{\upsilon}(\Ca)+ \Ca_{\upsilon}^{-1}\right)^{-1}\leq \phi_{t}(Q)\leq \Oa_\upsilon^{-1}+\Ca_{\upsilon}(\Oa).
$$
In addition, this implies
$$
\left(\Oa_{\upsilon}(\Ca)+ \Ca_{\upsilon}^{-1}\right)^{-1}\leq P\leq \Oa_\upsilon^{-1}+\Ca_{\upsilon}(\Oa)
	\quad\mbox{and}\quad
\left(\Oa_\upsilon^{-1}+\Ca_{\upsilon}(\Oa)\right)^{-1}\leq P^{-1}\leq\Oa_{\upsilon}(\Ca)+ \Ca_{\upsilon}^{-1}.
$$
\end{theo}

The following corollary is immediate.
\begin{cor}
For any $Q\in \SB_{r_1}^+$ and for any $t\geq \upsilon$ we have
$$
\mbox{\rm Spec}(\phi_{t}(Q))\quad\mbox{and}\quad \mbox{\rm Spec}(P)\quad
\subset\quad \left[ \left(\varpi_+^o(\Ca)+1/\varpi_-^c\right)^{-1},\varpi_+^c(\Oa)+1/\varpi_-^o\right]
$$
\end{cor}

Note that this corollary, together with the definition \eqref{ref-phi-uniform-max}, and following the proof of Proposition \ref{prop-lipschitz-continuity}, yields the following growth estimate,
$$
	\sup_{t\geq\upsilon}~t^{-1}\Vert \phi_t(Q)-Q\Vert_2<\infty
$$
in the general setting with time-varying signal models.

The proof of the theorem is based on comparison inequalities between the Riccati flow and the flow of matrices defined below. 

We let 
$$
Q\mapsto \phi^{o}_t(Q)\quad\mbox{\rm and}\quad Q\mapsto \phi^{c}_t(Q)
$$
with the flows associated with the Riccati equation with drift functions
$
\ricc^o$ and  $\ricc^c$ defined by
\begin{eqnarray*}
\ricc^c(Q)&:=&A_t Q+ QA_t^{\prime}+R_1\\
\ricc^o(Q)&:=&A_t Q+ QA_t^{\prime}-QS_tQ=\ricc(Q)-R_1.
\end{eqnarray*}

\begin{lem}
For any $t\geq  \upsilon$ we have
\begin{equation}\label{explicit-phico}
\Ca_t\leq \phi^{c}_t(Q)=\Ea_t(A)~Q~\Ea_t(A)^\prime+ \Ca_t
\quad\mbox{and}\quad   
\phi^{o}_t(Q)=\Ea_t(A)~\left(Q^{-1}+\overline{\Oa}_t\right)^{-1}~\Ea_t(A)^\prime\leq \Oa_t^{-1}
\end{equation}
with
$$
\overline{\Oa}_t:= \Ea_t(A)^\prime\Oa_t~\Ea_t(A) =\int_0^t \Ea_{s}(A)^\prime S_s~\Ea_{s}(A)~ds.
$$
In addition, for any $t\geq  \upsilon$ we have the estimates
$$
\phi^{o}_t(Q)\leq  \phi_t(Q)\leq  \Oa_t^{-1}+\Ca_t(\Oa)
$$
as well as
$$
\sup_{t\geq \upsilon}\phi_{t}(Q)\leq \Oa_\upsilon^{-1}+\Ca_{\upsilon}(\Oa) \quad\mbox{and}\quad 0<\left( \Oa_{\upsilon}^{-1}+\Ca_{\upsilon}(\Oa)\right)^{-1}\leq \inf_{t\geq \upsilon} \phi_{t}(Q)^{-1}.
$$
\end{lem}
\proof 
The l.h.s. inequality of (\ref{explicit-phico}) is immediate. We check the r.h.s. inequality of (\ref{explicit-phico}) using the fact that
\begin{eqnarray*}
	\partial_t\phi^{o}_t(Q) &=& \left(\partial_t \Ea_t(A)\right)\left(Q^{-1}+\overline{\Oa}_t\right)^{-1}\Ea_t(A)^\prime\\
			&&\hskip2cm+\Ea_t(A)\left(Q^{-1}+\overline{\Oa}_t\right)^{-1}\left(\partial_t\Ea_t(A)^\prime\right)\\
			&&\hskip4cm+\Ea_t(A)\left[\partial_t\left(Q^{-1}+\overline{\Oa}_t\right)^{-1}\right]\Ea_t(A)^\prime \\
		&=&A_t\phi^{o}_t(Q)+\phi^{o}_t(Q)A_t^{\prime}+{\Ea_t(A)\left[\partial_t\left(Q^{-1}+\overline{\Oa}_t\right)^{-1}\right]\Ea_t(A)^\prime}.
\end{eqnarray*}
On the other hand, recalling the inverse derivation formula
$$
\partial_tM^{-1}_t=-M^{-1}_t~\left(\partial_tM_t\right)~M^{-1}_t
$$
we find via Leibniz's rule that
$$
\begin{array}{l}
\partial_t\left(Q^{-1}+\overline{\Oa}_t\right)^{-1}\\
\\
\quad\quad=
-\left(Q^{-1}+\overline{\Oa}_t\right)^{-1}\underbrace{\left[\partial_t\left(Q^{-1}+\overline{\Oa}_t\right)
\right]}_{=\Ea_t(A)^\prime S_t~\Ea_t(A)}\left(Q^{-1}+\overline{\Oa}_t\right)^{-1}\\
\\
\quad\quad=-\left\{\left(Q^{-1}+\overline{\Oa}_t\right)^{-1}\Ea_t(A)^\prime\right\}S_t
\left\{\Ea_t(A)\left(Q^{-1}+\overline{\Oa}_t\right)^{-1}\right\}.
\end{array}
$$
This implies that
$$
\begin{array}{l}
{\Ea_t(A)\left[\partial_t\left(Q^{-1}+\overline{\Oa}_t\right)^{-1}\right]\Ea_t(A)^\prime}\\
\\
\quad =-{\left\{\Ea_t(A)\left(Q^{-1}+\overline{\Oa}_t\right)^{-1}\Ea_t(A)^\prime\right\}}~S_t~
{\left\{\Ea_t(A)\left(Q^{-1}+\overline{\Oa}_t\right)^{-1}\Ea_t(A)^\prime\right\}}
=-{\phi^{o}_t(Q)}S_t{\phi^{o}_t(Q)}.
\end{array}
$$
We also have
\begin{eqnarray*}
	\Ea_t(A)\left(Q^{-1}+\overline{\Oa}_t\right)^{-1}\Ea_t(A)^\prime &\leq&
		\Ea_t(A)\,\overline{\Oa}_t^{-1}\Ea_t(A)^\prime \\
		&=& \Ea_t(A)\,\Ea_t(A)^{-1}~\Oa_t^{-1}(\Ea_t(A)^{-1})^\prime\,\Ea_t(A)^\prime ~=~ \Oa_t^{-1}.
\end{eqnarray*}
This ends the proof of (\ref{explicit-phico}).
Also observe that
$$
\ricc(Q_1)-\ricc^o(Q_2)=\ricc(Q_1)-\ricc(Q_2)+R_1.
$$
Using the polarization formulae (\ref{polarization-formulae}) we conclude that
$$
\phi_t(Q)-\phi_t^o(Q)
\displaystyle=\int_0^t~\Ea_{s,t}(M(Q))~R_1~\Ea_{s,t}(M(Q))^{\prime}~ds\geq 0
$$
with the flow of matrices
$$
u\mapsto M_u(Q):=A_u-\frac{\phi^{o}_u(Q)+\phi_u(Q)}{2}S_u.
$$
We have the decomposition
\begin{eqnarray*}
\overline{\phi}_{t}(Q)&:=&\displaystyle \Ea_{t}(A)^{-1}\phi_{t}(Q)\Ea^\prime_{t}(A)^{-1}\\
	&=&\displaystyle Q+\int_0^{t} \Ea_{s}(A)^{-1}R_1 \Ea^\prime_{s}(A)^{-1}~ds \\
	&& \quad\quad-\int_0^{t}\Ea_{s}(A)^{-1} \phi_{s}(Q)\Ea^\prime_{s}(A)^{-1} \partial_s\overline{\Oa}_s~\Ea_{s}(A)^{-1} \phi_{s}(Q)\Ea^\prime_{s}(A)^{-1}~ds.
\end{eqnarray*}
In differential form this equation resumes to 
$$
\displaystyle\partial_t\overline{\phi}_{t}(Q) =\overline{R}_t -\overline{\phi}_{t}(Q)\left[\partial_t\overline{\Oa}_t\right]\overline{\phi}_{t}(Q)
\quad\mbox{\rm with $\overline{R}_t=\Ea_{t}(A)^{-1}R_1\Ea^\prime_{t}(A)^{-1}$.}
$$
On the other hand, we have
$$
\begin{array}{l}
\partial_t\left\{\overline{\Oa}_t\overline{\phi}_{t}(Q)
\overline{\Oa}_t\right\}\\
\\
\quad\quad=\left[\partial_t\overline{\Oa}_t\right] \overline{\phi}_{t}(Q)
\overline{\Oa}_t
+\overline{\Oa}_t\overline{\phi}_{t}(Q)
\left[\partial_t\overline{\Oa}_t\right]
+\overline{\Oa}_t\left\{\overline{R}_t
-\overline{\phi}_{t}(Q)~\left[\partial_t\overline{\Oa}_t\right]~ \overline{\phi}_{t}(Q)\right\}\overline{\Oa}_t\\
\\
\quad\quad=\overline{\Oa}_t\overline{R}_t\overline{\Oa}_t+\partial_t\overline{\Oa}_t
-\left[Id-\overline{\Oa}_t\overline{\phi}_{t}(Q)\right]\partial_t\overline{\Oa}_t
\left[Id-\overline{\Oa}_t\overline{\phi}_{t}(Q)\right]^{\prime}\\
\\
\quad\quad\leq \partial_t\overline{\Oa}_t+ \overline{\Oa}_t\overline{R}_t\overline{\Oa}_t:=\partial_t\left[\overline{\Oa}_t+\overline{\Ra}_t\right]
\end{array}
$$
with
$$
\overline{\Ra}_t=\int_0^t~\overline{\Oa}_s\overline{R}_s\overline{\Oa}_s~ds=\int_0^t~\Ea_s(A)^\prime\Oa_sR_1~\Oa_s~\Ea_s(A)~ds.
$$
This implies that
$$
\Ea_t(A)^\prime\Oa_t \phi_{t}(Q)\Oa_t \Ea_t(A)\,=\,\overline{\Oa}_t\overline{\phi}_{t}(Q)
\overline{\Oa}_t ~\leq \begin{array}[t]{rcl}\overline{\Oa}_t+\overline{\Ra}_t
&=\displaystyle\overline{\Oa}_t+\int_0^t\overline{\Oa}_s\overline{R}_s\overline{\Oa}_s~ds \qquad\qquad\qquad\qquad\qquad\\
&=\displaystyle \Ea_t(A)^\prime\Oa_t\Ea_t(A)+\int_0^t\Ea_s(A)^\prime\Oa_sR_1~\Oa_s\Ea_s(A)ds
\end{array}
$$
from which we conclude that
$$
 \phi_{t}(Q)\leq \Oa_t^{-1}+\Ca_t(\Oa).
$$
Since $\Ca_t(\Oa)$ and $\Oa_t$ don't depend on the initial state $Q$, for any $t\geq \upsilon$ we have
$$
 \phi_{t}(Q)= \phi_{\upsilon}(\phi_{t-\upsilon}(Q))\leq \Oa_\upsilon^{-1}+\Ca_{\upsilon}(\Oa).
$$
The inverse of both sides exists due to our observability/controllability assumption. This ends the proof of the lemma.
\qed

Whenever it exists the inverse $\phi_t(Q)^{-1}$ of the positive definite symmetric matrices
$
\phi_t(Q)>0
$, satisfies the following eigenvalue relationships
$$
\phi_t(Q)^{-1}\geq\left( \Oa_{\upsilon}^{-1}+\Ca_{\upsilon}(\Oa)\right)^{-1}\quad\Longrightarrow \inf_{Q\in\SB^+_{r_1}}\lambda_{min}\left(\phi_t(Q)^{-1}\right)\geq \left(\varpi_+^c(\Oa)+1/\varpi_-^o\right)^{-1}>0
$$
and 
$$
\sup_{Q\in\SB^+_{r_1} }\varsigma(\phi_t(Q))\leq \varsigma(\Oa_{\upsilon}^{-1}+\Ca_{\upsilon}(\Oa))\leq  \varsigma(\Oa_{\upsilon}^{-1})+\varsigma(\Ca_{\upsilon}(\Oa)) \leq \varpi^c_+(\Oa)+(1/\varpi^{o}_-).
$$ 
We also have
$$
\inf_{Q\in\SB^+_{r_1}}{\lambda_{min}\left(\phi_t(Q)^{-1}\right)}=\frac{1}{\sup_{Q\in\SB^+_{r_1}}{\lambda_{max}\left(\phi_t(Q)\right)}}\Rightarrow \sup_{Q\in\SB^+_{r_1}}{\lambda_{max}\left(\phi_t(Q)\right)}\leq \varpi_+^c(\Oa)+1/\varpi^o_-<\infty.
$$
The inverse  matrices $\phi_{t}(Q)^{-1}$ satisfy the equation
$$
\partial_t\phi_{t}(Q)^{-1}=\ricc_{-}(\phi_{t}(Q)^{-1})
$$
for any $t\geq \upsilon$, with the drift function
$$
\ricc_{-}(Q)=-A_t^{\prime}Q-QA_t-QR_1Q+S_t.
$$
We denote by $\phi^{-o}_t(Q^{-1})$ the flow starting at $Q^{-1}$ associated with the drift function
$$
\ricc_{-o}(Q):=-QA_t-A_t^{\prime}Q-QR_1Q=\ricc_{-}(Q)-S_t.
$$
The next lemma concerns the uniform boundedness of the inverse $\phi_t(Q)^{-1}$ of the positive definite symmetric matrices
$
\phi_t(Q)>0
$, w.r.t the time horizon $t\geq \upsilon$. 

\begin{lem}
For any $t\geq \upsilon$ we have
$$
\Ca_t\leq  \left[\phi^{-o}_t(Q^{-1})\right]^{-1}=\Ea_t(A) Q^{-1}\Ea_t(A)^\prime+\Ca_t=\phi^{c}_t(Q^{-1}).
$$
In addition we have
$$
\phi^{-o}_t(Q^{-1})\leq \phi_t(Q)^{-1}\leq \Oa_t(\Ca)+ \Ca_t^{-1}
$$
as well as
$$
\sup_{t\geq \upsilon}\phi_{t}(Q)^{-1}\leq  \Oa_{\upsilon}(\Ca)+ \Ca_{\upsilon}^{-1} \quad\mbox{and}\quad 0<\left(\Oa_{\upsilon}(\Ca)+ \Ca_{\upsilon}^{-1}\right)^{-1}\leq \inf_{t\geq \upsilon}{\phi_{t}(Q)}.
$$
\end{lem}
\proof
For brevity we write $\Ea_{s,t}:=\Ea_{s,t}(A)$. Arguing as in the proof of the preceding lemma, the flow $\phi^{-o}_t(Q^{-1})$ associated with $\ricc_{-o}$ is given by
\begin{eqnarray*}
 \phi^{-o}_t(Q^{-1})&=&(\Ea_t^{\prime})^{-1}~\left(Q^{-1}+\int_0^t\Ea_s^{-1}R_1~(\Ea_s^\prime)^{-1}~ds \right)^{-1}~\Ea_t^{-1}\\
	&=&\left(\Ea_t Q^{-1}\Ea_t^\prime+\int_0^t \Ea_{s,t}R_1~\Ea_{s,t}^\prime~ds \right)^{-1}\\
	&=&\left(\Ea_t Q^{-1}\Ea_t^\prime+\Ca_t
\right)^{-1}=\phi^{c}_t(Q^{-1})^{-1}~\leq~ \Ca_t^{-1}.
\end{eqnarray*}
This can be checked via differentiation as in the preceding proof. In addition, arguing as before,
$$
\phi_t(Q)^{-1}\geq\phi^{-o}_t(Q^{-1}).
$$
Also observe that the Riccati flow $\widehat{\phi}_t(Q):=\phi_t(Q)^{-1}$ satisfies a Riccati equation
defined similarly to that of $\phi_t(Q)$ but with a replacement on the matrices $(A_t,R_1,S_t)$ given by $(\widehat{A}_t,\widehat{R}_t,\widehat{S})$ with
$$
\widehat{A}_t:=-A_t^{\prime}\qquad \widehat{R}_t=S_t\quad \mbox{\rm and}\quad  \widehat{S}=R_1.
$$
That is, we have that
$$
\widehat{\phi}_t(Q)=\widehat{A}_t~\widehat{\phi}_t(Q)+\widehat{\phi}_t(Q)~\widehat{A}_t^{\prime}+\widehat{R}_t-\widehat{\phi}_t(Q)~
\widehat{S}~\widehat{\phi}_t(Q)
$$
with the initial condition $\widehat{\phi}_0(Q)=Q^{-1}$. This follows the inverse derivation formula also used in the preceding proof. Now it follows from the preceding lemma that
$$
\widehat{\phi}_t(Q)\leq  \widehat{\Oa}_t^{-1}+\widehat{\Ca}_t(\widehat{\Oa})
$$
with
$$
\widehat{\Oa}_t=\int_0^t \Ea^\prime_{s}(\widehat{A})^{-1}~R_1~\Ea_{s}(\widehat{A})^{-1}~ds= \int_0^t \Ea_s~R_1~\Ea_s^{\prime}~ds=\Ca_t
$$
and
\begin{eqnarray*}
\widehat{\Ca}_t(\widehat{\Oa})&=&
\widehat{\Oa}_t^{-1}\left[\int_0^t~\Ea_{s,t}^\prime(\widehat{A})^{-1}\widehat{\Oa}_s \widehat{R}_t~\widehat{\Oa}_s~\Ea_{s,t}(\widehat{A})^{-1}~ds\right]\widehat{\Oa}_t^{-1}\\
&=&\Ca_t^{-1}\left[\int_0^t~\Ea_{s,t}~\Ca_s~ S_t~\Ca_s~\Ea_{s,t}^\prime~ds\right]\Ca_t^{-1}=\Oa_t(\Ca).
\end{eqnarray*}
We conclude that
$$
\phi^{-o}_t(Q^{-1})\leq \phi_t(Q)^{-1}\leq \Ca_t^{-1}+ \Oa_t(\Ca).
$$
We also have
$$
\sup_{t\geq 0}\phi_{t}(Q)^{-1}=\sup_{t\geq \upsilon}\phi_{\upsilon}\left(\phi_{t-\upsilon}(Q)\right)^{-1}\leq \Oa_{\upsilon}(\Ca)+ \Ca_{\upsilon}^{-1}
$$
and therefore
$$
\inf_{t\geq \upsilon}{\phi_{t}(Q)}\geq\left(\Oa_{\upsilon}(\Ca)+ \Ca_{\upsilon}^{-1}\right)^{-1}.
$$
This ends the proof of the lemma.\qed

Combining this pair of lemmas we readily prove Theorem~\ref{compact-localisation}. This ends the proof of the theorem.

\subsection{Bucy's convergence theorem for Kalman-Bucy semigroups}~\label{bucy-theo-section}

We now prove the exponential convergence of a time-varying semigroup generated by a time-varying matrix difference of the form $A_t-\phi_t(Q)S_t$. This significantly generalises (\ref{hypothesis}) of Lemma \ref{cor-hypothesis}. Note here we seek explicit constants in this theorem which rely on our assumption $0<R_1\in\SB_{r_1}^+$. This assumption ensures controllability trivially holds, and it is a typical assumption in filtering applications. Qualitative exponential convergence estimates under weaker (i.e. only controllability/observability conditions on the model $(A_t,R_1,S_t)$) follow from a time-varying version of Lyapunov's method \cite{anderson67}; see e.g. \cite{bucy2,anderson67,anderson71} and \cite[Chapter 10]{Bittanti91}. The latter is of interest when $R_1$ may be only positive semi-definite, which is common in control applications, and more generally when considering the convergence properties of the Riccati equation for its own sake.

\begin{theo}[Bucy~\cite{bucy2}]\label{theo-expo-sg-ricc}
For any $t\geq s\geq \upsilon$ we have
$$
\sup_{Q\in\SB^+_{r_1}}\left\Vert  E_{s,t}(Q)\right\Vert_2
\leq \alpha\exp{\left\{-\beta (t-s)\right\}}
$$
with the parameters
$$
\alpha^2:= \frac{\varpi_+^o(\Ca)+1/\varpi_-^c}{\varpi_+^c(\Oa)+1/\varpi_-^o}\quad\mbox{and}\quad
2\beta :=\frac{1}{(\varpi_+^o(\Ca)+1/\varpi_-^c)}
 \left[\inf_{t\geq 0}\lambda_{min}(S_t)+\frac{\lambda_{min}(R_1)}{(\varpi_+^c(\Oa)+1/\varpi_-^o)^2}\right].
$$
\end{theo}

\proof
Observe that
$$
\begin{array}{l}
\phi_{t}(Q)^{-1}(A_t-\phi_{t}(Q)S_t)+(A_t-\phi_{t}(Q)S_t)^{\prime}\phi_{t}(Q)^{-1}+\ricc_{-}\left(\phi_{t}(Q)^{-1}\right) \qquad\qquad\qquad\qquad
	\\ \qquad\qquad\qquad\qquad\qquad\qquad\qquad\qquad\qquad\qquad\qquad\qquad\, =\, -\left[S_t+\phi_{t}(Q)^{-1}R_1\phi_{t}(Q)^{-1}\right].
\end{array}
$$
This implies that
$$
\begin{array}{l}
\partial_t \left(E_{s,t}(Q)^{\prime}\phi_{t}(Q)^{-1} E_{s,t}(Q)\right)\\
\\
\quad =E_{s,t}(Q)^{\prime}\left\{(A_t-\phi_{t}(Q)S_t)^{\prime}\phi_{t}(Q)^{-1} +\phi_{t}(Q)^{-1}(A_t-\phi_{t}(Q)S_t) +\ricc_{-}\left(
\phi_{t}(Q)^{-1}\right)\right\}E_{s,t}(Q)\end{array}
$$
from which we conclude that
$$
\partial_t \left(E_{s,t}(Q)^{\prime}\phi_{t}(Q)^{-1} E_{s,t}(Q)\right)=-E_{s,t}(Q)^{\prime}
\left[S_t+\phi_{t}(Q)^{-1}R_1\phi_{t}(Q)^{-1}\right]E_{s,t}(Q).
$$
By Theorem~\ref{compact-localisation}, we also have
\begin{eqnarray*}
	S_t+\phi_{t}(Q)^{-1}R_1~\phi_{t}(Q)^{-1} &\geq& \left[\inf_{t\geq 0}\lambda_{min}(S_t)+\frac{\lambda_{min}(R_1)}{\lambda^2_{max}(\phi_{t}(Q))} \right]~Id \\
 	&\geq&  \left[\inf_{t\geq 0}\lambda_{min}(S_t)+\frac{\lambda_{min}(R_1)}{(\varpi_+^c(\Oa)+1/\varpi_-^o)^2} \right]~Id
\end{eqnarray*}
and
$$
\begin{array}{l}
(\varpi_+^o(\Ca)+1/\varpi_-^c)^{-1}~Id~\leq \phi_{t}(Q)\leq (\varpi_+^c(\Oa)+1/\varpi_-^o)~Id~\\
\\
\qquad\quad\Longleftrightarrow\quad
(\varpi_+^o(\Ca)+1/\varpi_-^c)~Id~\geq \phi_{t}(Q)^{-1} \geq  (\varpi_+^c(\Oa)+1/\varpi_-^o)^{-1}~Id.
\end{array}$$
This implies that
$$
S_t+\phi_{t}(Q)^{-1}R_1~\phi_{t}(Q)^{-1}
\geq  \beta~\phi_{t}(Q)^{-1}
 $$
 with
 $$
 \beta:=\frac{1}{(\varpi_+^o(\Ca)+1/\varpi_-^c)}
 \left[\inf_{t\geq 0}\lambda_{min}(S_t)+\frac{\lambda_{min}(R_1)}{(\varpi_+^c(\Oa)+1/\varpi_-^o)^2}\right]
 $$
 from which we conclude that
 $$
 E_{s,t}(Q)^{\prime}\phi_{t}(Q)^{-1} E_{s,t}(Q)~\geq~ (\varpi_+^c(\Oa)+1/\varpi_-^o)~ E_{s,t}(Q)^{\prime}E_{s,t}(Q).
 $$
This implies that
$$
\partial_t\langle  E_{s,t}(Q)~x,\phi_{t}(Q)^{-1} E_{s,t}(Q)~x\rangle
~\leq~ -\beta~\langle
E_{s,t}(Q)~x, \phi_{t}(Q)^{-1}~E_{s,t}(Q)~x\rangle.
$$
By Gr\" onwall inequality we prove that
\begin{eqnarray*}
(\varpi_+^c(\Oa)+1/\varpi_-^o) \langle E_{s,t}(Q)~x,E_{s,t}(Q)~x\rangle
	&\leq& \langle E_{s,t}(Q)~x,\phi_{t}(Q)^{-1} E_{s,t}(Q)~x\rangle\\
	&\leq& e^{-\beta (t-s)} \langle x, \phi_{s}(Q)^{-1}~x\rangle\\
	&\leq& (\varpi_+^o(\Ca)+1/\varpi_-^c)~e^{-\beta (t-s)} \langle x, x\rangle
\end{eqnarray*}from which we conclude that
$$
	\Vert E_{s,t}(Q)\Vert_2^2\leq \frac{\varpi_+^o(\Ca)+1/\varpi_-^c}{\varpi_+^c(\Oa)+1/\varpi_-^o}~e^{-\beta (t-s)}.
$$
This ends the proof of the theorem.
\qed

We also have the following corollary. 
\begin{cor}\label{cor-bucy-any-t}
For any $0\leq s\leq t$ and any $Q\in\SB^+_{r_1}$ we have
\begin{equation}\label{final-ref-Est}
	\left\Vert  E_{s,t}(Q)\right\Vert_2 \leq  \rho(Q)\, \exp{\left\{-\beta (t-s)\right\}}
\end{equation}
with the function
\begin{eqnarray*}
Q\mapsto \rho(Q) :=( \alpha \vee 1)~\exp{\left[\left(\beta +\sup_{t\geq 0}\Vert A_t\Vert_2+\Vert \phi(Q)\Vert_2~\sup_{t\geq 0}\Vert S_t\Vert_2 \right)\upsilon \right]}
\end{eqnarray*}
and the uniform norm $\Vert \phi(Q)\Vert_2$  introduced in (\ref{ref-phi-uniform-max}).
\end{cor}
\proof
The estimate is immediate when $\upsilon\leq s\leq t$.
By Lemma~\ref{perturbation-lemma-intro}, for any  $0\leq s\leq t\leq \upsilon$ we have
$$
\left\Vert  E_{s,t}(Q)\right\Vert_2\leq \exp{\left\{-\beta (t-s)\right\}}~\exp{\left[\beta (t-s)+\left(\sup_{t\geq 0}\Vert A_t\Vert_2 +\Vert \phi(Q)\Vert_2~ \sup_{t\geq 0}\Vert S_t\Vert_2 \right)\upsilon\right]}~\Longrightarrow~ (\ref{final-ref-Est}).
$$
In the same vein, when $0\leq s\leq \upsilon\leq t$ we use theorem~\ref{theo-expo-sg-ricc} to check that 
\begin{eqnarray*}
E_{s,t}(Q)=E_{\upsilon,t}(Q)E_{s,\upsilon}(Q)\qquad\qquad\qquad\qquad\qquad\qquad\qquad\qquad\qquad\qquad\qquad\qquad\qquad\qquad \\
	\qquad\Longrightarrow~\Vert E_{s,t}(Q)\Vert_2\leq 
 \alpha\exp{\left\{-\beta (t-\upsilon)\right\}}~\exp{\left[\left(\sup_{t\geq 0}\Vert A_t\Vert_2 +\Vert \phi(Q)\Vert_2~ \sup_{t\geq 0}\Vert S_t\Vert_2 \right)(\upsilon-s)\right]}.
\end{eqnarray*}
This implies that
$$
\Vert E_{s,t}(Q)\Vert_2\leq 
 \alpha\exp{\left\{-\beta (t-s)\right\}}~\exp{\left[\left(\beta+\sup_{t\geq 0}\Vert A_t\Vert_2 +\Vert \phi(Q)\Vert_2~ \sup_{t\geq 0}\Vert S_t\Vert_2 \right)(\upsilon-s)\right]}.
 $$
This ends the proof of the corollary.
\qed

Using (\ref{polarization-formulae-phi1}) we readily check the following contraction estimate.
\begin{cor}\label{cor-lip-contraction-rsg}
For any $t\geq 0$ any $Q_1,Q_2\in\SB^+_{r_1}$ we have
$$
	\Vert\phi_t(Q_1)-\phi_t(Q_2)\Vert_2 \leq  \rho(Q_1,Q_2)\, \exp{\left\{-2\beta t\right\}}~\Vert Q_1-Q_2\Vert_2
$$
with 
$$
 \rho(Q_1,Q_2):= \rho(Q_1)\rho(Q_2)
$$
with $Q\mapsto \rho(Q)$ defined in Corollary~\ref{cor-bucy-any-t}.
\end{cor}

\subsection{Quantitative contraction estimates for time-invariant signal models}

We now consider time-invariant signal models, and we note that satisfaction of the observability and controllability rank conditions (\ref{def-contr-obs}) are sufficient to ensure the existence of a (unique) positive definite, and stabilizing, solution $P$ of (\ref{steady-state-eq}). We thus assume that the time-invariant matrix $A-PS$ satisfies (\ref{hypothesis}) of Lemma \ref{cor-hypothesis} for some $ \nu>0$~and some $\kappa<\infty$.

Bucy's theorem discussed in Section~\ref{bucy-theo-section} yields more or less directly several contraction inequalities. Note that because of (\ref{ref-phi-maxQ}) it follows that
\begin{eqnarray}
	Q\mapsto \rho(Q) &:=& (\alpha\vee 1)\,\exp{\left[\left(\beta + \Vert A\Vert_2+\Vert \phi(Q)\Vert_2~\Vert S\Vert_2\right)\,\upsilon \right]} \nonumber\\
		&\leq & (\alpha\vee 1)\,\exp{\left[\left(\beta + \Vert A\Vert_2+\left( \Vert P\Vert_2+\kappa^2\Vert Q-P\Vert_2\right)~\Vert S\Vert_2\right)\,\upsilon \right]} \label{cor-bucy-any-t-rho-homo}
\end{eqnarray}
and thus Corollaries \ref{cor-bucy-any-t} and \ref{cor-lip-contraction-rsg} immediately deliver crude (in terms of the rate) quantitative contraction results in the time-invariant setting. 

Note that the explicit constants following in this work rely again on our (common in filtering) assumption that $0<R_1\in\SB_{r_1}^+$. As before, qualitative convergence results, however of the same-nature, follow immediately; see \cite{anderson79,Bittanti91}. The following explicit estimates have now also been refined in \cite{bd-CARE}, without the assumption $0<R_1\in\SB_{r_1}^+$, which is of interest in more general applications of the Riccati equation.

The first result of this section concerns the exponential rate of the convergence of the Riccati flow towards the steady state.

\begin{cor}\label{cor-lip-contraction-rsg-mu}
For any $t\geq0$, and any $Q\in\SB^+_{r_1}$ we have
\begin{eqnarray*}
\Vert\phi_t(Q)-P\Vert_2
&\leq& ~\kappa_{\phi}(Q)~
e^{-2\nu t}~\Vert Q-P\Vert_2
\end{eqnarray*}
with the parameters
\begin{eqnarray*}
\kappa_{\phi}(Q):=\kappa^2
~\exp{\left\{(2\beta )^{-1}~\Vert S\Vert_2~\kappa^2~ \rho(P,Q)~\Vert Q-P\Vert_2 ~\right\}}.
\end{eqnarray*}
In the above, $\rho(P,Q)$ is defined in Corollary~\ref{cor-bucy-any-t} and Corollary~\ref{cor-lip-contraction-rsg} noting the upper-bound \eqref{cor-bucy-any-t-rho-homo}.
\end{cor}
\proof
Combining (\ref{polarization-formulae-phi3}) and (\ref{hypothesis}) with Corollary~\ref{cor-bucy-any-t} and Corollary~\ref{cor-lip-contraction-rsg} we prove that
\begin{eqnarray*}
Z_t&:=&~e^{2\nu t}~\Vert\phi_t(Q)-P\Vert_2~\\
&\leq&~ \kappa^2
Z_0+ \kappa^2~\Vert S\Vert_2~ \rho(P,Q)~\Vert Q-P\Vert_2~\int_0^t~
~\exp{\left\{-2\beta u\right\}}~
Z_u~du
\end{eqnarray*}
for any $\epsilon\in ]0,1]$ and any $t\geq0$. 
A direct application of Gr\" onwall inequality yields 
$$
Z_t\leq ~\kappa^2~
Z_0~\exp{\left\{(2\beta )^{-1}\kappa^2~\Vert S\Vert_2~ \rho(P,Q)~\Vert Q-P\Vert_2\int_0^t~2\beta ~\exp{\left\{-2\beta u\right\}}~du\right\}}.
$$
This ends the proof of the corollary.
\qed

We also have the following quantitative exponential convergence result on the Kalman-Bucy semigroup and a contraction-type inequality on the Riccati flow.

\begin{cor}\label{cor-Est-varsigma}
For any $0\leq s\leq t$, $\epsilon\in ]0,1]$ and any $Q\in\SB^+_{r_1}$ we have
\begin{eqnarray*}
\Vert E_{s,t}(Q)\Vert_2
&\leq& \kappa_{E}(Q)\,\exp{\left[-\nu (t-s)\right]}~~\Longrightarrow~~ \Vert E(Q)\Vert_2~:=~
\sup_{0\leq s\leq t}\Vert E_{s,t}(Q)\Vert_2 ~\leq ~\kappa_{E}(Q)
\end{eqnarray*}
with the parameters
$$
	\kappa_{E}(Q):=\kappa~\exp{\left(\frac{\kappa}{2\nu}~  \kappa_{\phi}(Q)~\Vert S\Vert_2~\Vert Q-P\Vert_2~\right)}
$$
where $\kappa_{\phi}(Q)$ is the parameter defined in Corollary~\ref{cor-lip-contraction-rsg-mu}.

In addition, for any $Q_1,Q_2\in\SB^+_{r_1}$ we have the local Lipschitz property
\begin{equation}\label{phi-exact-varepsilon}
\Vert\phi_t(Q_1)-\phi_t(Q_2)\Vert_2\leq ~ \kappa_{\phi}(Q_1,Q_2)
\exp{\left[-2\nu t\right]}~\Vert Q_1-Q_2\Vert_2
\end{equation}
with
$$
 \kappa_{\phi}(Q_1,Q_2):=
 \kappa_{E}(Q_1) \kappa_{E}(Q_2).
$$
\end{cor}
\proof
For any $0\leq s\leq t$ we have
$$
	E_{s,t}(Q)=\exp{\left[\oint_s^t\left[(A-PS)+(P-\phi_u(Q))S\right]~du\right]}.
$$
By Corollary~\ref{cor-lip-contraction-rsg-mu} we have
\begin{eqnarray*}
\int_s^t\Vert (\phi_u(Q)-P)S\Vert_2\,du&\leq& \kappa_{\phi}(Q)~\Vert S\Vert_2~\Vert Q-P\Vert_2~
 \int_s^t
e^{-2\nu u}~du~\leq~ \frac{ \kappa_{\phi}(Q)~\Vert S\Vert_2}{2\nu}~\Vert Q-P\Vert_2.
\end{eqnarray*}
The first estimate is a direct consequence of Lemma~\ref{perturbation-lemma-intro}.
The estimate (\ref{phi-exact-varepsilon}) is a direct consequence of (\ref{polarization-formulae-phi1}).
The end of the proof is now a direct consequence of Lemma~\ref{perturbation-lemma-intro}.
\qed

The preceding result can be contrast with Corollaries \ref{cor-bucy-any-t} and \ref{cor-lip-contraction-rsg} that concern the more general time-varying Kalman-Bucy semigroup and Riccati flow (noting the upper-bound given by \eqref{cor-bucy-any-t-rho-homo}). In the time-invariant domain one may be able to improve the exponential rate significantly, e.g. via $\nu>0$ in (\ref{hypothesis}) of Lemma \ref{cor-hypothesis}, as compared to Bucy's general Theorem \ref{theo-expo-sg-ricc}. For example, $\nu$ may be related to the spectral abscissa of $A-PS$. However, any improvement in the rate may come at the expense of a (much) larger constant.

Tracking the nature of constants is important when applying these results in practice; e.g. when analysing the stability of ensemble Kalman filters \cite{tong16,dm-16-enkf}, or extended Kalman filters \cite{reif2000,dm-16-eekf}. These constants are also typically related to the underlying state-space dimension. In this sense, one should carefully follow the form and the source of these terms, in order to understand accurately the dimensional error dependence of any approximation (in, e.g., high-dimensions) \cite{mf-dm-04,rebeschini2015,tong-16-enkf}. As noted, the constants in the preceding rely on our assumption $0<R_1\in\SB_{r_1}^+$ which is common in filtering applications. In the time-invariant setting, these constants have been explicitly refined in \cite{bd-CARE} and without this assumption on $R_1$.

Finally, we have a quantitative contraction inequality on on the Kalman-Bucy semigroup with time-invariant signal models.

\begin{cor}\label{lipschitz-E}
For any $0\leq s\leq t$, $\epsilon\in ]0,1]$ and any $Q_1,Q_2\in\SB_{r_1}^+$ we have
$$
\Vert E_{s,t}(Q_2)-E_{s,t}(Q_1)\Vert_2\leq \kappa_E(Q_1,Q_2)
\exp{\left[-\nu (t-s)\right]}~\Vert Q_2-Q_1\Vert_2
$$
with
$$
\kappa_E(Q_1,Q_2):=\kappa_{E}(Q_2)~+\kappa_{E}^2(Q_2)~\kappa_{\phi}(Q_1,Q_2)~\Vert S\Vert_2~(2\nu)^{-1}.
$$
\end{cor}

\proof
We have
\begin{eqnarray*}
\partial_t\left(E_{s,t}(Q_2)-E_{s,t}(Q_1)\right)&=&(A-\phi_t(Q_2)S)~E_{s,t}(Q_2)-
(A-\phi_t(Q_1)S)~E_{s,t}(Q_1)\\
&=&(A-\phi_t(Q_2)S)~(E_{s,t}(Q_2)-E_{s,t}(Q_1))\\
&&\hskip4cm+(\phi_t(Q_1)-\phi_t(Q_2))SE_{s,t}(Q_1).
\end{eqnarray*}
This implies that
$$
E_{s,t}(Q_2)-E_{s,t}(Q_1)=E_{s,t}(Q_2)~(Q_2-Q_1)+\int_s^t E_{u,t}(Q_2)
(\phi_u(Q_1)-\phi_u(Q_2))SE_{s,u}(Q_1)~du.
$$
By Corollary~\ref{cor-Est-varsigma}, this yields the estimate
$$
\begin{array}{l}
\exp{\left[\nu (t-s)\right]}~\Vert E_{t}(Q_2)-E_{t}(Q_1)\Vert_2\\
\\
\qquad\qquad\displaystyle\leq \kappa_{E}(Q_2)~\Vert Q_2-Q_1\Vert_2+\kappa^2_{E}(Q_2)\Vert S\Vert_2\int_s^t ~
\Vert\phi_u(Q_1)-\phi_u(Q_2)\Vert_2 ~
du.
\end{array}
$$
Using (\ref{phi-exact-varepsilon}) we check that
$$
\begin{array}{l}
\exp{\left[\nu (t-s)\right]}~\Vert E_{s,t}(Q_2)-E_{s,t}(Q_1)\Vert_2\\
\\
\qquad\qquad\displaystyle\leq \Vert Q_2-Q_1\Vert_2\left[\kappa_{E}(Q_2)~+\kappa^2_{E}(Q_2)~ \kappa_{\phi}(Q_1,Q_2)~\Vert S\Vert_2~(2\nu)^{-1}\right].
\end{array}
$$
This ends the proof of the corollary.
\qed

\section{Contraction of Kalman-Bucy-type stochastic flows} \label{sec-diffusion-flows}

Firstly, we review a straightforward qualitative stability result for the time-varying Kalman-Bucy filter that follows from the uniform boundedness of the Riccati flow. 

\begin{theo} \label{conv-timevarying-KF}
For any $t\geq s\geq \upsilon$ we have the uniform estimate,
$$
\sup_{Q\in\SB^+_{r_1}}\left\Vert\,  \EE\left( \psi_{s,t}(x,Q) - \varphi_{s,t}(X_s) \,\vert\, X_s\right)\, \right\Vert_2
\leq \alpha\exp{\left\{-\beta (t-s)\right\}}\Vert\,x-X_s \Vert_2
$$
with the parameters $\alpha,\beta >0$ defined in Theorem \ref{theo-expo-sg-ricc}.

\end{theo}
\proof
The proof (and result) follows that of Theorem \ref{theo-expo-sg-ricc}. By uniform observability/controllability we suppose $Q=\phi_{s-\upsilon,s}(Q_0)$ without loss of generality. Let $\xi_{s,t} := \EE\left( \psi_{s,t}(x,Q) - \varphi_{s,t}(X_s) \,\vert\, X_s\right)$. Consider the functional
$$
	\frac{\left\Vert \xi_{s,t} \right\Vert_2^2}{\left(\varpi_+^c(\Oa)+1/\varpi_-^o\right)} ~\leq~ \xi_{s,t}' \,\phi_{s,t}(Q)^{-1}\xi_{s,t}
 ~\leq ~(\varpi_+^o(\Ca)+1/\varpi_-^c)\left\Vert  \xi_{s,t}\right\Vert_2^2.
$$
Then
\begin{eqnarray*}
	\partial_t~\xi_{s,t}\,'\phi_{s,t}(Q)^{-1}\xi_{s,t} &=& - \,\xi_{s,t}'\,\left(\phi_{s,t}(Q)^{-1}R_1\phi_{s,t}(Q)^{-1}+S_t\right) \xi_{s,t} \\
		&\leq& -\,\beta\, \xi_{s,t}\,'\phi_{s,t}(Q)^{-1}\xi_{s,t}.
\end{eqnarray*}
By Gr\" onwall's inequality we find
\begin{eqnarray*}
(\varpi_+^c(\Oa)+1/\varpi_-^o) \xi_{s,t}\,' \xi_{s,t} 
	&\leq& \xi_{s,t}\,'\phi_{s,t}(Q)^{-1}\xi_{s,t} \\
	&\leq& e^{-\beta (t-s)} \xi_{s,s}\,'\phi_{s,s}(Q)^{-1}\xi_{s,s} \\
	&\leq& (\varpi_+^o(\Ca)+1/\varpi_-^c)~e^{-\beta (t-s)} \xi_{s,s}\,'\xi_{s,s}
\end{eqnarray*}
and the result follows with $\alpha,\beta >0$ defined as in Theorem \ref{theo-expo-sg-ricc}.
\qed

Given this classical, qualitative, stability result, we now study in more precise terms the convergence of Kalman-Bucy stochastic flows, both in the classical filtering form, and the novel nonlinear diffusion form. We study exponential inequalities that bound, with dedicated probability, the stochastic flow of the sample paths at any time, with respect to the underlying signal. We also provide almost sure contraction-type estimates. Both types of result offer a notion of filter stability and the analysis in this section is novel. 

We assume that the signal models are time-invariant throughout the remainder of this section, and we build on the quantitative estimates of the prior section.

In further development of this section, going forward we consider the function
$$
Q\mapsto \sigma(Q):=2\sqrt{2}~\kappa_E(Q)~\left[\left(\Vert \phi(Q)\Vert_2^2~\Vert S\Vert_2+\Vert R\Vert_2\right)~r_1/\nu\right]^{1/2}
$$
as well as
 \begin{eqnarray*}
  \chi_0(Q_1,Q_2)&=&\nu^{-1}\kappa_{E}(Q_1) \kappa_{\phi}(Q_1,Q_2)\\
 \chi_1(Q)&=&\Vert S\Vert_2~\kappa_E(Q)/2\quad\mbox{and}\quad
 \chi_2(Q)=\Vert S\Vert_2~\sigma(Q)+2~\sqrt{2r_1\Vert S\Vert_2\nu}
 \end{eqnarray*}
with the parameters $\nu$, $\Vert \phi(Q)\Vert_2$, 
$\kappa_{E}(Q)$, and  $\kappa_{\phi}(Q_1,Q_2)$ respectively 
defined in the exponential rate of (\ref{hypothesis}), in 
Proposition~\ref{prop-lipschitz-continuity}, in Corollary~\ref{cor-Est-varsigma}, and in Corollary~\ref{cor-Est-varsigma} again.

\subsection{Time-invariant Kalman-Bucy filter}

\begin{theo}\label{prop-event-control-KB}
The conditional probability of the following
events
\begin{equation}\label{event-control-KB}
\left\Vert~\psi_{s,t}(x,Q)-\varphi_{s,t}(X_s)-E_{t-s}(Q)\left[x-X_s\right]~\right\Vert_2~
\leq\frac{e^2}{\sqrt{2}}~\left[\frac{1}{2}+\left(\delta+\sqrt{\delta}\right)\right]\,\sigma^2(Q)
\end{equation}
given the state variable $X_s$ is greater than $1-e^{-\delta}$, for any $\delta\geq 0$ and any $t\in[s,\infty[$.
\end{theo}

By (\ref{event-control-KB}), the conditional probability of the following
event
\begin{equation}\label{event-control-KB-bis}
\begin{array}[t]{l}
\left\Vert~\psi_{s,t}(x,Q)-\varphi_{s,t}(X_s)\right\Vert_2~\\
\\
\qquad \leq  \displaystyle\kappa_E(Q)~e^{-\nu(t-s)}~\Vert x-X_s\Vert_2+ \frac{e^2}{\sqrt{2}}~\left[\frac{1}{2}+\left(\delta+\sqrt{\delta}\right)\right]~\sigma^2(Q)
\end{array}
\end{equation}
given the state variable $X_s$ is greater than $1-e^{-\delta}$, for any $\delta\geq 0$ and any $t\in[s,\infty[$.

The above theorem is a direct consequence of (\ref{event-control}) and the following technical lemma.

\begin{lem}\label{lem-tex-ucontrol}
For any  $x\in \RR^{r_1}$, $Q\in\SB^+_{r_1}$ and $n\geq 1$ we have
the uniform estimate
$$
\sup_{t\geq s}\EE\left(\left\Vert~\psi_{s,t}(x,Q)-\varphi_{s,t}(X_s)-E_{t\vert s}(Q)~\left[x-X_s\right]~\right\Vert_2^{2n}~\vert~X_s\right)^{1/n}
\displaystyle\leq n~\sigma^2(Q).
$$
In particular, for any $t\geq s\geq0$. we have
$$
\EE\left(\left\Vert~\psi_{s,t}(x,Q)-\varphi_{s,t}(X_s)\right\Vert_2^{2n}~\vert~X_s\right)^{\frac{1}{2n}}
\displaystyle\leq \sqrt{n}~\sigma(Q)+\kappa_E(Q)~e^{-\nu(t-s)}~
\left\Vert x-X_s\right\Vert_2.
$$
\end{lem}
\proof
For any given $s\geq0$ and for any $t\in [s,\infty[$, and any $x\in \RR^{r_1}$ we have
\begin{eqnarray*}
d\left[\psi_{s,t}(x,Q)-\varphi_{s,t}(X_s)\right]&=&\left[A-\phi_{s,t}(Q)S\right]~\left[\psi_{s,t}(x,Q)-\varphi_{s,t}(X_s)\right]~dt+
dM_{s,t}
\end{eqnarray*}
with the $r_1$-multivariate martingale $(M_{s,t})_{t\in [s,\infty[}$ given by
\begin{equation}\label{def-Mst}
\begin{array}{l}
\displaystyle t\in [s,\infty[\mapsto M_{s,t}~=~\int_s^t\phi_{s,u}(Q) C^{\prime}R^{-1}_{2}dV_u-R_1^{1/2}~(W_t-W_s)~\\
\\
\qquad\qquad\qquad\qquad\Longrightarrow~~\partial_t\,\langle M_{s,\bullet}(k),M_{s,\bullet}(l)\rangle_t~=~\phi_{s,t}(Q)S\phi_{s,t}(Q)+R_1.
\end{array}\end{equation}
This yields the formula
$$
N_{s,t}:=\left[\psi_{s,t}(x,Q)-\varphi_{s,t}(X_s)\right]-E_{t\vert s}(Q)~\left[x-X_s\right]
=\int_s^t~E_{u-s,t-s}(Q)~
~dM_{s,u}.
$$

On the other hand, we have
$$
\begin{array}{l}
\displaystyle\EE\left(\Vert N_{s,t}
\Vert_2^{2n}\right)^{1/n}=\left[\EE\left(\left[\sum_{1\leq k\leq r_1} N_{s,t}(k)^2
\right]^n\right)\right]^{1/n} 
\displaystyle\leq \sum_{1\leq k\leq r_1} \EE\left(N_{s,t}(k)^{2n}\right)^{1/n} \\
\\
\qquad\qquad\displaystyle =\sum_{1\leq k\leq r_1} \EE\left(\left[\sum_{1\leq l\leq r_1}\int_s^t~E_{u-s,t-s}(Q)(k,l)
~dM_{s,u}(l)\right]^{2n}\right)^{1/n}.
\end{array}
$$

By the Burkolder-Davis-Gundy inequality (\ref{BDG}) we have
$$
\begin{array}{l}
\displaystyle\EE\left(\left[\sum_{1\leq l\leq r_1}
\int_s^t~E_{u,t\vert s}(Q)(k,l)~dM_{s,u}(l)\right]^{2n}\right)^{1/n}\\
\\ 
\qquad\qquad\displaystyle\leq 4^2n
\sum_{1\leq l,l^{\prime}\leq r_1}\int_s^t~E_{u,t\vert s}(Q)(k,l)E_{u,t\vert s}(Q)(k,l^{\prime})~(\phi_{s,{u}}(Q)S\phi_{s,{u}}(Q)+R_1)(l,l^{\prime})~du \\
\\
\qquad\qquad\displaystyle\leq 4^2n
\int_s^t~\left[E_{u,t\vert s}(Q)(\phi_{s,{u}}(Q)S\phi_{s,{u}}(Q)+R_1)E_{u,t\vert s}(Q)^{\prime}\right](k,k)~du.
\end{array}
$$
This implies that
\begin{eqnarray*}
\EE\left(\Vert N_{s,t}\Vert_2^{2n}\right)^{1/n}&\leq& 
4^2n
\int_s^t~\tr\left[(\phi_{s,{u}}(Q)S\phi_{s,{u}}(Q)+R_1)E_{u,t\vert s}(Q)^{\prime}E_{u,t\vert s}(Q)\right]~du\\
&\leq & 4^2n~\int_s^t~\Vert \phi_{s,{u}}(Q)S\phi_{s,{u}}(Q)+R_1\Vert_2~\Vert E_{u,t\vert s}(Q)\Vert_F^2~du.
\end{eqnarray*}
This yields
$$
\begin{array}{l}
\displaystyle\EE\left(\Vert \int_s^t~E_{u,t\vert s}(Q)~
~dM_{s,u}
\Vert_2^{2n}\right)^{1/n}
\displaystyle\leq 
4^2 n~\left(\Vert \phi(Q)\Vert_2^2~\Vert S\Vert_2+\Vert R_1\Vert_2\right)~r_1~\int_s^t~\Vert
E_{u-s,t-s}(Q)\Vert_2^2~du.
\end{array}
$$
This implies that
$$
\begin{array}{l}
\EE\left(\left\Vert\psi_{s,t}(x,Q)-\varphi_{s,t}(X_s)-E_{t\vert s}(Q)~\left[x-X_s\right]\right\Vert_2^{2n}\right)^{1/n}\\
\\
\qquad\qquad\displaystyle\leq 
4^2 n~\left(\Vert \phi(Q)\Vert_2^2~\Vert S\Vert_2+\Vert R_1\Vert_2\right)~r_1~\int_s^t~\Vert
E_{u-s,t-s}(Q)\Vert_2^2~du.\end{array}
$$
Using Proposition~\ref{prop-lipschitz-continuity} and Corollary~\ref{cor-Est-varsigma} we have

$$
\begin{array}{l}
\EE\left(\left\Vert\psi_{s,t}(x,Q)-\varphi_{s,t}(X_s)-E_{t\vert s}(Q)~\left[x-X_s\right]\right\Vert_2^{2n}\right)^{1/n}\\
\\
\qquad\qquad\displaystyle\leq 8 n\nu^{-1}\kappa_E^2(Q)~\left(\Vert \phi(Q)\Vert_2^2~\Vert S\Vert_2+\Vert R_1\Vert_2\right)~r_1~\int_s^t~(2\nu)~\exp{\left[-2\nu (t-u)\right]}~du.
\end{array}
$$
This ends the proof of the lemma.\qed

\begin{theo}\label{theo-psi-st}
For any $t\geq s\geq 0$, $x_1,x_2\in \RR^{r_1}$, $Q_1,Q_2\in\SB^+_{r_1}$ and $n\geq 1$ we have the almost sure local contraction estimate
 $$
\begin{array}{l}
\EE\left(\Vert\psi_{s,t}(x_1,Q_1)-\psi_{s,t}(x_2,Q_2)\Vert_2^{2n}~\vert~X_s\right)^{\frac{1}{2n}}
\\
\\
\qquad\qquad\qquad\leq \displaystyle \kappa_{E}(Q_1)~e^{-\nu (t-s)}~\Vert x_1-x_2 \Vert_2\\
\\
\qquad\qquad\qquad\qquad\displaystyle+e^{-\nu (t-s)}~  \chi_0(Q_1,Q_2)
\left\{ \chi_1(Q_2)~\left\Vert x_2-X_s\right\Vert_2~+~\sqrt{n}~ \chi_2(Q_1)\right\}~\Vert Q_1-Q_2\Vert_2.
\end{array}$$

\end{theo}

\proof
We have
$$
\begin{array}{l}
d\left(\psi_{s,t}(x_1,Q_1)-\psi_{s,t}(x_2,Q_2)\right)\\
\\
\qquad\qquad=\displaystyle\left\{\left[A-\phi_{s,t}(Q_1)S\right]~\psi_{s,t}(x_1,Q_1)~-~\left[A-\phi_{s,t}(Q_2)S\right]~\psi_{s,t}(x_2,Q_2)\right\}
dt\\
\\
\hskip8cm\displaystyle+\left[\phi_{s,t}(Q_1)-\phi_{s,t}(Q_2)\right]
~C^{\prime}R^{-1}_{2}~dY_t\\
\\
\qquad\qquad=\displaystyle\left[A-\phi_{s,t}(Q_1)S\right]~\left(\psi_{s,t}(x_1,Q_1)-\psi_{s,t}(x_2,Q_2)\right)~dt\\
\\
\displaystyle\hskip4cm-\left[\phi_{s,t}(Q_1)-\phi_{s,t}(Q_2)\right]S~\psi_{s,t}(x_2,Q_2)~dt\\
\\
\hskip6cm\displaystyle+\left[\phi_{s,t}(Q_1)-\phi_{s,t}(Q_2)\right]
~(S \varphi_{s,t}(X_s)~ds+ C^{\prime}R^{-1/2}_{2}~dV_t).
\end{array}
$$
This yields the decomposition
$$
\begin{array}{l}
d\left[\psi_{s,t}(x_1,Q_1)-\psi_{s,t}(x_2,Q_2)\right]\\
\\
\qquad\qquad=\displaystyle\left[A-\phi_{s,t}(Q_1)S\right]~\left[\psi_{s,t}(x_1,Q_1)-\psi_{s,t}(x_2,Q_2)\right]~dt\\
\\
\displaystyle\hskip4cm+\left[\phi_{s,t}(Q_1)-\phi_{s,t}(Q_2)\right]S~\left[\varphi_{s,t}(X_s)-\psi_{s,t}(x_2,Q_2)\right]~dt\\
\\
\hskip6cm\displaystyle+\left[\phi_{s,t}(Q_1)-\phi_{s,t}(Q_2)\right]~dM_t
\end{array}
$$
with $M_t=C^{\prime}R^{-1/2}_{2}~V_t$.
This yields the decomposition
$$
\begin{array}{l}
\psi_{s,t}(x_1,Q_1)-\psi_{s,t}(x_2,Q_2)\\
\\
\qquad=\displaystyle E_{t\vert s}(Q_1)~(x_1-x_2)+\int_s^t~E_{u,t\vert s}(Q_1)~
\left[\phi_{s,u}(Q_1)-\phi_{s,u}(Q_2)\right]S~\left[\varphi_{s,u}(X_s)-\psi_{s,u}(x_2,Q_2)\right]~du\\
\\
\hskip6.5cm\displaystyle+\int_s^t~E_{u,t\vert s}(Q_1)~\left[\phi_{s,u}(Q_1)-\phi_{s,u}(Q_2)\right]~dM_u.
\end{array}$$
Arguing as in the proof of Lemma~\ref{lem-tex-ucontrol} we have
$$
\begin{array}{l}
\displaystyle\EE\left(\Vert \int_s^t~E_{u,t\vert s}(Q_1)~\left[\phi_{s,u}(Q_1)-\phi_{s,u}(Q_2)\right]~dM_u
\Vert_2^{2n}\right)^{1/n}\\
\\
\qquad\qquad\qquad\displaystyle\leq  4^2 nr_1~\Vert S\Vert_2~\int_s^t~\Vert \phi_{s,u}(Q_1)-\phi_{s,u}(Q_2)\Vert_2^2
~\Vert E_{u,t\vert s}(Q)\Vert_2^2~du.
\end{array}
$$
 By Corollary~\ref{cor-Est-varsigma} we have the estimate
$$
\begin{array}{l}
\displaystyle\EE\left(\Vert \int_s^t~E_{u,t\vert s}(Q_1)~\left[\phi_{s,u}(Q_1)-\phi_{s,u}(Q_2)\right]~dM_u
\Vert_2^{2n}\right)^{1/n}\\
\\
\qquad\displaystyle\leq 
 4^2 nr_1~\Vert S\Vert_2~
 ~\kappa_{E}^2(Q_1)~ \kappa_{\phi}^2(Q_1,Q_2)~\Vert Q_1-Q_2\Vert^2_2~\exp{\left[-2\nu (t-s)\right]} ~\int_s^t~
~\exp{\left[-2\nu (u-s)\right]}
~du\\
\\
\qquad\displaystyle\leq 
8nr_1~(\Vert S\Vert_2/\nu)~
 ~(\kappa_{E}(Q_1)~ \kappa_{\phi}(Q_1,Q_2)~\Vert Q_1-Q_2\Vert_2)^2~\exp{\left[-2\nu (t-s)\right]} \end{array}
$$ 
from which we find that
$$
\begin{array}{l}
\displaystyle\EE\left(\Vert \int_s^t~E_{u,t\vert s}(Q_1)~\left[\phi_{s,u}(Q_1)-\phi_{s,u}(Q_2)\right]~dM_u
\Vert_2^{2n}\right)^{\frac{1}{2n}}\\
\\
\qquad\qquad\displaystyle\leq \sqrt{n}~
2\sqrt{2r_1(\Vert S\Vert_2/\nu)}~
 ~\kappa_{E}(Q_1)~ \kappa_{\phi}(Q_1,Q_2)~\Vert Q_1-Q_2\Vert_2~\exp{\left[-\nu (t-s)\right]}. \end{array}
$$ 
On the other hand we have the almost sure estimate
$$
\begin{array}{l}
\Vert\psi_{s,t}(x_1,Q_1)-\psi_{s,t}(x_2,Q_2)\Vert_2\\
\\
\qquad\leq \displaystyle \Vert E_{t\vert s}(Q_1)\Vert_2~\Vert x_1-x_2 \Vert_2\\
\\
\hskip1cm\displaystyle+\Vert S\Vert_2 
\int_s^t~\Vert E_{u,t\vert s}(Q_1)\Vert_2~
\Vert \phi_{s,u}(Q_1)-\phi_{s,u}(Q_2)\Vert_2~\Vert\varphi_{s,u}(X_s)-\psi_{s,u}(x_2,Q_2)\Vert_2~du\\
\\
\hskip6.5cm\displaystyle+\Vert \int_s^t~E_{u,t\vert s}(Q_1)~\left[\phi_{s,u}(Q_1)-\phi_{s,u}(Q_2)\right]~dM_u\Vert_2.
\end{array}$$
Combining Corollary~\ref{cor-Est-varsigma} and the estimate
 (\ref{phi-exact-varepsilon}) 
 we prove that
 $$
\begin{array}{l}
\Vert\psi_{s,t}(x_1,Q_1)-\psi_{s,t}(x_2,Q_2)\Vert_2\\
\\
\qquad\leq \displaystyle \kappa_{E}(Q_1)~\exp{\left[-\nu (t-s)\right]}~\Vert x_1-x_2 \Vert_2+\Vert S\Vert_2~\kappa_{E}(Q_1)~ \kappa_{\phi}(Q_1,Q_2)~\Vert Q_1-Q_2\Vert_2~\\
\\
\qquad\qquad\displaystyle\times
\int_s^t~~\exp{\left[-\nu (t-u)\right]}
\exp{\left[-2\nu (u-s)\right]}~\Vert\varphi_{s,u}(X_s)-\psi_{s,u}(x_2,Q_2)\Vert_2~du\\
\\
\hskip5cm\displaystyle+\Vert \int_s^t~E_{u-s,t-s}(Q_1)~\left[\phi_{s,u}(Q_1)-\phi_{s,u}(Q_2)\right]~dM_u\Vert_2.
\end{array}$$
 This implies that
  $$
\begin{array}{l}
\Vert\psi_{s,t}(x_1,Q_1)-\psi_{s,t}(x_2,Q_2)\Vert_2\\
\\
\qquad\leq \displaystyle \kappa_{E}(Q_1)~\exp{\left[-\nu (t-s)\right]}~\Vert x_1-x_2 \Vert_2+\Vert S\Vert_2~\kappa_{E}(Q_1)~ \kappa_{\phi}(Q_1,Q_2)~\Vert Q_1-Q_2\Vert_2~\\
\\
\qquad\qquad\displaystyle\times~~\exp{\left[-\nu (t-s)\right]}~
\int_s^t~\exp{\left[-\nu (u-s)\right]}~\Vert\varphi_{s,u}(X_s)-\psi_{s,u}(x_2,Q_2)\Vert_2~du\\
\\
\hskip5cm\displaystyle+\Vert \int_s^t~E_{u-s,t-s}(Q_1)~\left[\phi_{s,u}(Q_1)-\phi_{s,u}(Q_2)\right]~dM_u\Vert_2.
\end{array}$$
Using the generalized Minkowski inequality we check that
  $$
\begin{array}{l}
\EE\left(\Vert\psi_{s,t}(x_1,Q_1)-\psi_{s,t}(x_2,Q_2)\Vert_2^{2n}~\vert~X_s\right)^{\frac{1}{2n}}\\
\\
\qquad\leq \displaystyle \kappa_{E}(Q_1)~\exp{\left[-\nu (t-s)\right]}~\Vert x_1-x_2 \Vert_2+\Vert S\Vert_2~\kappa_{E}(Q_1)~ \kappa_{\phi}(Q_1,Q_2)~\Vert Q_1-Q_2\Vert_2~\\
\\
\qquad\qquad\displaystyle\times~~\exp{\left[-\nu (t-s)\right]}~
\int_s^t~\exp{\left[-\nu (u-s)\right]}~\EE\left(\Vert\varphi_{s,u}(X_s)-\psi_{s,u}(x_2,Q_2)\Vert_2^{2n}~\vert~X_s\right)^{\frac{1}{2n}}~du\\
\\
\hskip5cm\displaystyle+\EE\left(\Vert \int_s^t~E_{u-s,t-s}(Q_1)~\left[\phi_{s,u}(Q_1)-\phi_{s,u}(Q_2)\right]~dM_u\Vert_2^{2n}\right)^{\frac{1}{2n}}.
\end{array}$$
 Lemma~\ref{lem-tex-ucontrol} combined with the Burkholder-Davis-Gundy estimates stated above 
 implies that
  $$
\begin{array}{l}
\exp{\left[\nu (t-s)\right]}~\EE\left(\Vert\psi_{s,t}(x_1,Q_1)-\psi_{s,t}(x_2,Q_2)\Vert_2^{2n}~\vert~X_s\right)^{\frac{1}{2n}}\\
\\
\qquad\leq \displaystyle \kappa_{E}(Q_1)~\Vert x_1-x_2 \Vert_2+\Vert S\Vert_2~\kappa_{E}(Q_1)~ \kappa_{\phi}(Q_1,Q_2)~\Vert Q_1-Q_2\Vert_2~\\
\\
\qquad\qquad\displaystyle\times~
\int_s^t~\exp{\left[-\nu (u-s)\right]}~\left\{ \sqrt{n}~\sigma(Q_2)+\kappa_E(Q_2)~e^{-\nu(u-s)}~
\left\Vert x_2-X_s\right\Vert_2\right\}~du\\
\\
\hskip5cm\displaystyle+ \sqrt{n}~
2\sqrt{2r_1(\Vert S\Vert_2/\nu)}~
 ~\kappa_{E}(Q_1)~ \kappa_{\phi}(Q_1,Q_2)~\Vert Q_1-Q_2\Vert_2.
\end{array}$$
 Observe that
   $$
\begin{array}{l}
\displaystyle \int_s^t~\exp{\left[-\nu (u-s)\right]}~\left\{ \sqrt{n}~\sigma(Q_2)+\kappa_E(Q_2)~e^{-\nu(u-s)}~
\left\Vert x_2-X_s\right\Vert_2\right\}~du\\
\\
\qquad\qquad\qquad= \sqrt{n}~(\sigma(Q_2)/\nu)
+(\kappa_E(Q_2)/(2\nu))~\left\Vert x_2-X_s\right\Vert_2.
 \end{array}
 $$
This yields
   $$
\begin{array}{l}
\exp{\left[\nu (t-s)\right]}~\EE\left(\Vert\psi_{s,t}(x_1,Q_1)-\psi_{s,t}(x_2,Q_2)\Vert_2^{2n}~\vert~X_s\right)^{\frac{1}{2n}}\\
\\
\qquad\leq \displaystyle \kappa_{E}(Q_1)~\Vert x_1-x_2 \Vert_2+\Vert S\Vert_2~\kappa_{E}(Q_1)~ \kappa_{\phi}(Q_1,Q_2)~\Vert Q_1-Q_2\Vert_2~\\
\\
\qquad\qquad\displaystyle\times\left[\sqrt{n}~(\sigma(Q_2)/\nu)
+(\kappa_E(Q_2)/(2\nu))~\left\Vert x_2-X_s\right\Vert_2\right]\\
\\
\hskip5cm\displaystyle+ \sqrt{n}~
2\sqrt{2r_1(\Vert S\Vert_2/\nu)}~
 ~\kappa_{E}(Q_1)~ \kappa_{\phi}(Q_1,Q_2)~\Vert Q_1-Q_2\Vert_2
\end{array}$$
from which we conclude that
  $$
\begin{array}{l}
\exp{\left[\nu (t-s)\right]}~\EE\left(\Vert\psi_{s,t}(x_1,Q_1)-\psi_{s,t}(x_2,Q_2)\Vert_2^{2n}~\vert~X_s\right)^{\frac{1}{2n}}
\leq \displaystyle \kappa_{E}(Q_1)~\Vert x_1-x_2 \Vert_2\\
\\
\qquad\qquad\qquad\qquad\displaystyle+
\left\{\left[ 
\Vert S\Vert_2~\kappa_E(Q_2)/2~\right]~\left\Vert x_2-X_s\right\Vert_2~+~\sqrt{n}~\left[\Vert S\Vert_2~\sigma(Q_2)+2\sqrt{2r_1\Vert S\Vert_2\nu}~
 \right]\right\}\\
 \\
 \qquad\qquad\qquad\qquad\qquad\qquad\displaystyle\times~\nu^{-1}\kappa_{E}(Q_1) \kappa_{\phi}(Q_1,Q_2)~\Vert Q_1-Q_2\Vert_2.
\end{array}$$
This ends the proof of the theorem.\qed

\subsection{Nonlinear time-invariant Kalman-Bucy diffusions}

\begin{theo}\label{prop-event-control-KB-diffusion}

For any $t\geq s\geq 0$, $x\in \RR^{r_1}$, $Q\in\SB^+_{r_1}$ and $n\geq 1$ we have
$$
\EE\left(\left\Vert~\overline{\psi}_{s,t}(x,Q)-\varphi_{s,t}(X_s)\right\Vert_2^{2n}~\vert~X_s\right)^{\frac{1}{2n}}
\displaystyle\leq \sqrt{2n}~\sigma(Q)+\kappa_E(Q)~e^{-\nu(t-s)}~
\left\Vert x-X_s\right\Vert_2.
$$
The conditional probability of the following
events
\begin{equation}
\left\Vert~\overline{\psi}_{s,t}(x,Q)-\varphi_{s,t}(X_s)-E_{t-s}(Q)~\left[x-X_s\right]~\right\Vert_2~
\leq \sqrt{2}~e^2~\left[\frac{1}{2}+\left(\delta+\sqrt{\delta}\right)\right]~\sigma^2(Q) \nonumber
\end{equation}
given the state variable $X_s$ is greater than $1-e^{-\delta}$, for any $\delta\geq 0$ and any $t\in[s,\infty[$.
\end{theo}

\proof
Observe that
\begin{eqnarray*}
d\left[
\overline{\psi}_{s,t}(x,Q)-\varphi_{s,t}(X_s)\right]&=&\left[A-\phi_{s,t}(Q)S\right]~
\left[\overline{\psi}_{s,t}(x,Q)-\varphi_{s,t}(X_s)\right]~dt+d\overline{M}_{s,t}
\end{eqnarray*}
with an $r_1$-valued martingale $(\overline{M}_{s,t})_{t\geq s}$ defined by
$$
t\in [s,\infty[\mapsto \begin{array}[t]{rcl}
\overline{M}_{s,t}&=&\displaystyle\int_s^t~\phi_{s,u}(Q)~C^{\prime}R^{-1/2}_{2}d(V_u-\overline{V}_u)+R^{1/2}_1\left[(\overline{W}_t-W_t)-(\overline{W}_s-W_s)\right]\\
&\stackrel{law}{=}&\sqrt{2}~M_{s,t}
\end{array}
$$ 
with the martingale $(M_{s,t})_{t\in [s,\infty[}$ discussed in (\ref{def-Mst}). The proof now follows the same arguments as the 
proof of Theorem~\ref{prop-event-control-KB} and Lemma~\ref{lem-tex-ucontrol}, thus it is skipped.
\qed

In the same vein, recalling that
\begin{eqnarray*}
d\overline{\psi}_{s,t}(x,Q)&=&\left[\left[A-\phi_{s,t}(Q)S\right]~\overline{\psi}_{s,t}(x,Q)+\phi_{s,t}(Q)~S~\varphi_{s,t}(X_s)\right]~dt\\
&&\qquad\qquad\qquad\qquad+R^{1/2}_{1}~d\overline{W}_t+\phi_{s,t}(Q)C^{\prime}R^{-1/2}_{2}d(V_t-\overline{V}_{t})
\end{eqnarray*}
we find the decomposition
$$
\begin{array}{l}
\overline{\psi}_{s,t}(x_1,Q_1)-\overline{\psi}_{s,t}(x_2,Q_2)\\
\\
\qquad=\displaystyle E_{t\vert s}(Q_1)~(x_1-x_2)+\int_s^t~E_{u,t\vert s}(Q_1)~
\left[\phi_{s,u}(Q_1)-\phi_{s,u}(Q_2)\right]S~\left[\varphi_{s,u}(X_s)-\overline{\psi}_{s,u}(x_2,Q_2)\right]~du\\
\\
\qquad\qquad\qquad\displaystyle+\int_s^t~E_{u,t\vert s}(Q_1)~\left[\phi_{s,u}(Q_1)-\overline{\psi}_{s,u}(Q_2)\right]~d\overline{M}_u
\end{array}$$
with $\overline{M}_t=\sqrt{2}~C^{\prime}R^{-1/2}_{2}~(V_t-\overline{V}_t)/\sqrt{2}$. 

\begin{theo} \label{theo-psi-diffusion-st}
For any $t\geq s\geq 0$, $x_1,x_2\in \RR^{r_1}$, $Q_1,Q_2\in\SB^+_{r_1}$ and $n\geq 1$ we have
the almost sure local contraction estimate
 $$
\begin{array}{l}
\EE\left(\Vert \overline{\psi}_{s,t}(x_1,Q_1)-\overline{\psi}_{s,t}(x_2,Q_2)\Vert_2^{2n}~\vert~X_s\right)^{\frac{1}{2n}}
\\
\\
\quad\quad\leq \displaystyle \kappa_{E}(Q_1)~e^{-\nu (t-s)}~\Vert x_1-x_2 \Vert_2\\
\\
\qquad\qquad\displaystyle+  \sqrt{2}~e^{-\nu (t-s)}~\chi_0(Q_1,Q_2)
\left\{ \chi_1(Q_2)~\left\Vert x_2-X_s\right\Vert_2~+~\sqrt{n}~\chi_2(Q_1)\right\}~\Vert Q_1-Q_2\Vert_2.
\end{array}$$
\end{theo}

\noindent Proof of this theorem follows readily that of Theorem~\ref{theo-psi-st}.

The analysis in this section encapsulates and extends the existing convergence and stability results for the Kalman-Bucy filter; e.g. as studied in \cite{kalman61,bucy68,anderson71,ocone-pardoux}. And we capture the properties of this convergence in a more quantitative manner than previously considered. The use of our nonlinear Kalman-Bucy diffusion provides a novel interpretation of the Kalman-Bucy filter that allows one to consider a more general class of signal model in a natural manner.

In particular, stability of the nonlinear Kalman-Bucy diffusion implies convergence of the filter, given arbitrary initial conditions, to the conditional mean of the signal given the observation filtration. Moreover, it implies convergence of the conditional distribution to a Gaussian defined by the conditional mean of the Kalman-Bucy diffusion and its covariance. Similar results were considered by Ocone and Pardoux in \cite{ocone-pardoux} but with no quantitative analysis. 

Note that our analysis further provides exponential relationships between the actual sample paths of the filter and the signal (with dedicated probability).

This analysis completes our review of the Kalman-Bucy filter and its stability properties.

\section*{Acknowledgements}
Adrian N. Bishop is with the University of Technology Sydney (UTS) and Data61 (CSIRO), Analytics Research Group. He is supported by the Australian Research Council (ARC) via a Discovery Early Career Researcher Award (DE-120102873). He is also an adjunct Fellow at the Australian National University (ANU).

Pierre Del Moral is with INRIA, Bordeaux Research Center (France) and the University of New South Wales (UNSW), School of Mathematics and Statistics (Australia).


\begin{thebibliography}{99.}


\bibitem{abou-kandil03}
H. Abou-Kandil, G. Freiling, V. Ionescu, and G. Jank. Matrix Riccati Equations in Control and Systems Theory. Birkhauser Verlag (2003).

\bibitem{anderson71}
B.D.O. Anderson. Stability properties of Kalman-Bucy filters. Journal of the Franklin Institute. vol. 291, no. 2. pp. 137--144 (1971).

\bibitem{anderson67}
B.D.O. Anderson and J.B. Moore. Time-Varying Version of the Lemma of Lyapunov. Electronics Letters. vol. 3, no. 7. pp. 293--294 (1967).

\bibitem{anderson79}
B.D.O. Anderson and J.B. Moore. Optimal Filtering. Dover Publications (1979).

\bibitem{bd-CARE}
A.N. Bishop and P. Del Moral. An explicit Floquet-type representation of Riccati aperiodic exponential semigroups. arXiv e-print, \href{https://arxiv.org/abs/1805.02127}{\tt arXiv:1805.02127} (2018).

\bibitem{Bittanti91}
S. Bittanti, A.J. Laub and J.C. Willems (Editors). The Riccati Equation. Springer-Verlag (1991).

\bibitem{bucy2}
R.S. Bucy. Global Theory of the Riccati Equation. Journal of Computer and System Sciences. vol. 1. pp. 349--361 (1967).

\bibitem{bucy72corrected}
R.S. Bucy. The Riccati Equation and Its Bounds. Journal of Computer and System Sciences. vol. 6. pp. 342--353 (1972).

\bibitem{bucy72remarks}
R.S. Bucy. Remarks on ``A Note on Bounds on Solutions of the Riccati Equation''. IEEE Transactions on Automatic Control. vol 17, no. 1. pp. 179 (1972).

\bibitem{bucy75}
R.S. Bucy. Structural Stability for the Riccati Equation. SIAM Journal on Control. vol. 13, no. 4. pp. 749--753 (1975).

\bibitem{bucy68}
R.S. Bucy and P.D. Joseph. Filtering for Stochastic Processes with Applications to Guidance. Interscience Publishers (1968).

\bibitem{bucy72}
R.S. Bucy and J. Rodriguez-Canabal. A Negative Definite Equilibrium and its Induced Cone of Global Existence for the Riccati Equation. SIAM Journal on Mathematical Analysis. vol 3, no.4. pp. 644--646 (1972).

\bibitem{callier81}
F.M. Callier and J.L. Willems. Criterion for the Convergence of the Solution of the Riccati Differential Equation. IEEE Transactions on Automatic Control. vol. 26, no. 6. pp. 1232--1242 (1981).

\bibitem{callier95}
F.M. Callier and J. Winkin. Convergence of the Time-Invariant Riccati Differential Equation towards Its Strong Solution for Stabilizable Systems. Journal of Mathematical Analysis and Applications. vol. 192, no. 1. pp. 230--257 (1995).

\bibitem{nicolao92}
G. De Nicolao and M. Gevers. Difference and Differential Riccati Equations: A Note on the Convergence to the Strong Solution. IEEE Transactions on Automatic Control. vol. 37, no. 7. pp. 1055--1057 (1992).

\bibitem{mf-dm-04}
P. Del Moral. Feynman-Kac Formulae. Springer (2004).

\bibitem{mf-dm-13}
P. Del Moral. Mean field simulation for Monte Carlo integration. Chapman \& Hall/CRC Monographs on Statistics \& Applied Probability (2013).

\bibitem{dm-16-eekf}
P. Del Moral, A. Kurtzmann, J. Tugaut. On the stability and the uniform propagation of chaos of extended ensemble Kalman-Bucy filters. arXiv e-print, arXiv:1606.08256 (2016).

\bibitem{dm-16-enkf}
P. Del Moral and J. Tugaut. On the stability and the uniform propagation of chaos properties of ensemble Kalman-Bucy filters. Annals of Applied Probability. vol. 28, no. 2. pp 790--850 (2018). arXiv e-print, \href{https://arxiv.org/pdf/1605.09329.pdf}{\tt arXiv:1605.09329}.

\bibitem{delyon2001}
B. Delyon. A note on uniform observability. IEEE Transactions on Automatic Control. vol. 46, no. 8. pp. 1326--1327 (2001).

\bibitem{evensen-review}
G. Evensen. The Ensemble Kalman Filter: theoretical formulation and practical implementation. Ocean Dynamics. vol. 53, pp. 343--367 (2003).

\bibitem{fielder}
M. Fiedler. Special Matrices and Their Applications in Numerical Mathematics. 2nd edition. Dover Publications (2008). 

\bibitem{gevers85}
M. Gevers, R.R. Bitmead, I.R. Petersen and R.J. Kaye. When is the solution of the Riccati equation stabilizing at every iteration?. In `Frequency domain and state space methods for linear systems' edited by C.I. Byrnes and A. Lindquist. pp. 531--540 (1986).

\bibitem{hitz72}
K.L. Hitz, T.E. Fortmann, and B.D.O. Anderson. A Note on Bounds on Solutions of the Riccati Equation. IEEE Transactions on Automatic Control. vol. 17, no. 1. pp. 178 (1972).

\bibitem{kalman60}
R.E. Kalman. Contributions to the Theory of Optimal Control. Boletin de la Sociedad Matematica Mexicana. vol. 5. pp. 102--119 (1960).

\bibitem{kalman60-2}
R. E. Kalman. A New Approach to Linear Filtering and Prediction Problems. Journal of Basic Engineering. vol. 82, no. 1. pp. 35--45 (1960).

\bibitem{kalman72}
R.E. Kalman. Further Remarks on ``A Note on Bounds on Solutions of the Riccati Equation''. IEEE Transactions on Automatic Control. vol 17, no. 1. pp. 179--180 (1972).

\bibitem{kalman60-3}
R.E. Kalman and J.E. Bertram. Control Systems Analysis and Design Via the ``Second Method'' of Lyapunov. Journal of Basic Engineering. vol. 82, no. 2. pp. 371--393 (1960).

\bibitem{kalman61}
R.E. Kalman and R.S. Bucy. New Results in Linear Filtering and Prediction Theory. Journal of Basic Engineering. vol. 83, no. 1. pp. 95--108 (1961).

\bibitem{Khoshnevisan}
D. Khoshnevisan. Analysis of Stochastic Partial Differential Equations. American Mathematical Soc. (2014).

\bibitem{kresemir}
V. Kresimir. Damped Oscillations of Linear Systems: A Mathematical Introduction. Springer-Verlag. Lecture Notes in Mathematics (2011).

\bibitem{kucera72}
V. Kucera. A Contribution to Matrix Quadratic Equations. IEEE Transactions on Automatic Control. vol. 17, no. 3. pp. 344--347 (1972).

\bibitem{kucera73}
V. Kucera. A Review of the Matrix Riccati Equation. Kybernetika. vol. 9, no. 1. pp. 42--61 (1973).

\bibitem{Kwakernaak72}
H. Kwakernaak and R. Sivan. Linear Optimal Control Systems. Wiley-Interscience (1972).

\bibitem{Lancaster1995}
P. Lancaster and L. Rodman. Algebraic Riccati Equations. Oxford University Press (1995).

\bibitem{legland09}
F. Le Gland, V. Monbet, V.-D. Tran. Large sample asymptotics for the ensemble Kalman filter. Research Report: RR-7014, INRIA. $<$inria-00409060$>$ (2009).

\bibitem{tong-16-enkf}
A.J. Majda and X.T. Tong. Robustness and accuracy of finite ensemble Kalman filters in large dimensions. arXiv e-print, arXiv:1606.09321 (2016).

\bibitem{vanloan-19}
C. Moler, C. Van Loan. Nineteen dubious ways to compute the exponential of a matrix.
SIAM Review. vol. 45, no. 1. pp. 3--49 (2003).

\bibitem{Molinari73}
B.P. Molinari. The stabilizing solution of the algebraic Riccati equation. SIAM Journal on Control. vol. 11, no. 2. pp. 262--271 (1973).

\bibitem{Molinari77}
B.P. Molinari. The time-invariant linear-quadratic optimal control problem. Automatica. vol. 13, no. 4. pp. 347--357 (1977).

\bibitem{ocone-pardoux}
D. Ocone, E. Pardoux.
Asymptotic stability of the optimal filter with respect to its initial condition. SIAM Journal on Control and Optimization. vol. 34, no. 1. pp. 226-243 (1996).

\bibitem{Park97}
P. Park and T. Kailath. Convergence of the DRE solution to the ARE strong solution. IEEE Transactions on Automatic Control. vol. 42, no. 4. 573--578 (1997). 

\bibitem{Poubelle88}
M-A. Poubelle, R.R. Bitmead and M.R. Gevers. Fake Algebraic Riccati Techniques and Stability. IEEE Transactions on Automatic Control. vol. 33, no. 4. pp. 379--381 (1988).

\bibitem{Poubelle86}
M-A. Poubelle, I.R. Petersen, M.R. Gevers, and R.R. Bitmead. A Miscellany of Results on an Equation of Count J. F. Riccati. IEEE Transactions on Automatic Control. vol. 31, no. 7. pp. 651--654 (1986).

\bibitem{rebeschini2015}
P. Rebeschini and R. Van Handel. Can local particle filters beat the curse of dimensionality?. The Annals of Applied Probability. vol. 25, no. 5. pp. 2809--2866 (2015).

\bibitem{reid72}
W.T. Reid. Riccati Differential Equations. Academic Press. (1972).

\bibitem{reif2000}
K. Reif, S. Gunther, E. Yaz, and R. Unbehauen. Stochastic stability of the continuous-time extended Kalman filter. IEE Proceedings - Control Theory and Applications. vol. 147, no. 1. pp. 45--52 (2000).

\bibitem{yao-feng-reng}
Y.F. Ren. On the Burkholder-Davis-Gundy inequalities for continuous martingales. Statistics and Probability Letters.
vol. 78. pp. 3034-3039 (2008).

\bibitem{shayman86}
M.A. Shayman. Phase Portrait of the Matrix Riccati Equation. SIAM Journal on Control and Optimization. vol. 24, no. 1. pp. 1--65 (1986).

\bibitem{tong16}
X.T. Tong, A.J. Majda and D. Kelly. Nonlinear stability and ergodicity of ensemble based Kalman filters. Nonlinearity. vol. 29, no. 2. pp. 657--691 (2016).

\bibitem{vanloan}
C. Van Loan. The sensitivity of the matrix exponential. SIAM Journal on Numerical Analysis. vol. 14, no. 6. pp. 971--981 (1977).

\bibitem{willems71}
J.C. Willems. Least Squares Stationary Optimal Control and the Algebraic Riccati Equation. IEEE Transactions on Automatic Control. vol. 16, no. 6. pp. 621--634 (1971).

\bibitem{wimmer85}
H.K Wimmer. Monotonicity of maximal solutions of algebraic Riccati equations. Systems \& Control Letters. vol. 5, no. 5. pp. 317--319 (1985).

\bibitem{wonham68}
W.M. Wonham. On a Matrix Riccati Equation of Stochastic Control. SIAM Journal of Control. vol. 6, no. 4. pp 681--697 (1968).


\end{thebibliography}
\end{document}